\magnification=\magstep1
\def\firstpage{1}
\pageno=\firstpage
\font\fiverm=cmr5
\font\sevenrm=cmr7
\font\sevenbf=cmbx7
\font\eightrm=cmr8
\font\eightbf=cmbx8
\font\ninerm=cmr9
\font\ninebf=cmbx9
\font\tenbf=cmbx10
\font\twelvebf=cmbx12

\font\magnineeufm=eufm9 scaled\magstep1

%
%
\newskip\ttglue
\font\fiverm=cmr5
\font\fivei=cmmi5
\font\fivesy=cmsy5
\font\fivebf=cmbx5
\font\sixrm=cmr6
\font\sixi=cmmi6
\font\sixsy=cmsy6
\font\sixbf=cmbx6
\font\sevenrm=cmr7
\font\eightrm=cmr8
\font\eighti=cmmi8
\font\eightsy=cmsy8
\font\eightit=cmti8
\font\eightsl=cmsl8
\font\eighttt=cmtt8
\font\eightbf=cmbx8
\font\ninerm=cmr9
\font\ninei=cmmi9
\font\ninesy=cmsy9
\font\nineit=cmti9
\font\ninesl=cmsl9
\font\ninett=cmtt9
\font\ninebf=cmbx9
\font\twelverm=cmr12
\font\twelvei=cmmi12
\font\twelvesy=cmsy12
\font\twelveit=cmti12
\font\twelvesl=cmsl12
\font\twelvett=cmtt12
\font\twelvebf=cmbx12


\def\eightpoint{\def\rm{\fam0\eightrm}
  \textfont0=\eightrm \scriptfont0=\sixrm \scriptscriptfont0=\fiverm
  \textfont1=\eighti  \scriptfont1=\sixi  \scriptscriptfont1=\fivei
  \textfont2=\eightsy  \scriptfont2=\sixsy  \scriptscriptfont2=\fivesy
  \textfont3=\tenex  \scriptfont3=\tenex  \scriptscriptfont3=\tenex
  \textfont\itfam=\eightit  \def\it{\fam\itfam\eightit}
  \textfont\slfam=\eightsl  \def\sl{\fam\slfam\eightsl}
  \textfont\ttfam=\eighttt  \def\tt{\fam\ttfam\eighttt}
  \textfont\bffam=\eightbf  \scriptfont\bffam=\sixbf
    \scriptscriptfont\bffam=\fivebf  \def\bf{\fam\bffam\eightbf}
  \tt  \ttglue=.5em plus.25em minus.15em
  \normalbaselineskip=9pt
  \setbox\strutbox=\hbox{\vrule height7pt depth2pt width0pt}
  \let\sc=\sixrm  \let\big=\eightbig \normalbaselines\rm}

\def\eightbig#1{{\hbox{$\textfont0=\ninerm\textfont2=\ninesy
        \left#1\vbox to6.5pt{}\right.$}}}


\def\ninepoint{\def\rm{\fam0\ninerm}
  \textfont0=\ninerm \scriptfont0=\sixrm \scriptscriptfont0=\fiverm
  \textfont1=\ninei  \scriptfont1=\sixi  \scriptscriptfont1=\fivei
  \textfont2=\ninesy  \scriptfont2=\sixsy  \scriptscriptfont2=\fivesy
  \textfont3=\tenex  \scriptfont3=\tenex  \scriptscriptfont3=\tenex
  \textfont\itfam=\nineit  \def\it{\fam\itfam\nineit}
  \textfont\slfam=\ninesl  \def\sl{\fam\slfam\ninesl}
  \textfont\ttfam=\ninett  \def\tt{\fam\ttfam\ninett}
  \textfont\bffam=\ninebf  \scriptfont\bffam=\sixbf
    \scriptscriptfont\bffam=\fivebf  \def\bf{\fam\bffam\ninebf}
  \tt  \ttglue=.5em plus.25em minus.15em
  \normalbaselineskip=11pt
  \setbox\strutbox=\hbox{\vrule height8pt depth3pt width0pt}
  \let\sc=\sevenrm  \let\big=\ninebig \normalbaselines\rm}

\def\ninebig#1{{\hbox{$\textfont0=\tenrm\textfont2=\tensy
        \left#1\vbox to7.25pt{}\right.$}}}


\def\twelvepoint{\def\rm{\fam0\twelverm}
  \textfont0=\twelverm \scriptfont0=\eightrm \scriptscriptfont0=\sixrm
  \textfont1=\twelvei  \scriptfont1=\eighti  \scriptscriptfont1=\sixi
  \textfont2=\twelvesy  \scriptfont2=\eightsy  \scriptscriptfont2=\sixsy
  \textfont3=\tenex  \scriptfont3=\tenex  \scriptscriptfont3=\tenex
  \textfont\itfam=\twelveit  \def\it{\fam\itfam\twelveit}
  \textfont\slfam=\twelvesl  \def\sl{\fam\slfam\twelvesl}
  \textfont\ttfam=\twelvett  \def\tt{\fam\ttfam\twelvett}
  \textfont\bffam=\twelvebf  \scriptfont\bffam=\eightbf
    \scriptscriptfont\bffam=\sixbf  \def\bf{\fam\bffam\twelvebf}
  \tt  \ttglue=.5em plus.25em minus.15em
  \normalbaselineskip=11pt
  \setbox\strutbox=\hbox{\vrule height8pt depth3pt width0pt}
  \let\sc=\sevenrm  \let\big=\twelvebig \normalbaselines\rm}

\def\twelvebig#1{{\hbox{$\textfont0=\tenrm\textfont2=\tensy
        \left#1\vbox to7.25pt{}\right.$}}}
\catcode`\@=11
%

\def\undefine#1{\let#1\undefined}
\def\newsymbol#1#2#3#4#5{\let\next@\relax
 \ifnum#2=\@ne\let\next@\msafam@\else
 \ifnum#2=\tw@\let\next@\msbfam@\fi\fi
 \mathchardef#1="#3\next@#4#5}
\def\mathhexbox@#1#2#3{\relax
 \ifmmode\mathpalette{}{\m@th\mathchar"#1#2#3}%
 \else\leavevmode\hbox{$\m@th\mathchar"#1#2#3$}\fi}
\def\hexnumber@#1{\ifcase#1 0\or 1\or 2\or 3\or 4\or 5\or 6\or 7\or 8\or
 9\or A\or B\or C\or D\or E\or F\fi}

\font\tenmsa=msam10
\font\sevenmsa=msam7
\font\fivemsa=msam5
\newfam\msafam
\textfont\msafam=\tenmsa
\scriptfont\msafam=\sevenmsa
\scriptscriptfont\msafam=\fivemsa
\edef\msafam@{\hexnumber@\msafam}
\mathchardef\dabar@"0\msafam@39
\def\dashrightarrow{\mathrel{\dabar@\dabar@\mathchar"0\msafam@4B}}
\def\dashleftarrow{\mathrel{\mathchar"0\msafam@4C\dabar@\dabar@}}

\def\ulcorner{\delimiter"4\msafam@70\msafam@70 }
\def\urcorner{\delimiter"5\msafam@71\msafam@71 }
\def\llcorner{\delimiter"4\msafam@78\msafam@78 }
\def\lrcorner{\delimiter"5\msafam@79\msafam@79 }
\def\yen{{\mathhexbox@\msafam@55}}
\def\checkmark{{\mathhexbox@\msafam@58}}
\def\circledR{{\mathhexbox@\msafam@72}}
\def\maltese{{\mathhexbox@\msafam@7A}}

\font\tenmsb=msbm10
\font\sevenmsb=msbm7
\font\fivemsb=msbm5
\newfam\msbfam
\textfont\msbfam=\tenmsb
\scriptfont\msbfam=\sevenmsb
\scriptscriptfont\msbfam=\fivemsb
\edef\msbfam@{\hexnumber@\msbfam}
\def\Bbb#1{{\fam\msbfam\relax#1}}
\def\widehat#1{\setbox\z@\hbox{$\m@th#1$}%
 \ifdim\wd\z@>\tw@ em\mathaccent"0\msbfam@5B{#1}%
 \else\mathaccent"0362{#1}\fi}
\def\widetilde#1{\setbox\z@\hbox{$\m@th#1$}%
 \ifdim\wd\z@>\tw@ em\mathaccent"0\msbfam@5D{#1}%
 \else\mathaccent"0365{#1}\fi}
\font\teneufm=eufm10
\font\seveneufm=eufm7
\font\fiveeufm=eufm5
\newfam\eufmfam
\textfont\eufmfam=\teneufm
\scriptfont\eufmfam=\seveneufm
\scriptscriptfont\eufmfam=\fiveeufm

\catcode`\@=11
\newsymbol\boxdot 1200
\newsymbol\boxplus 1201
\newsymbol\boxtimes 1202
\newsymbol\square 1003
\newsymbol\blacksquare 1004
\newsymbol\centerdot 1205
\newsymbol\lozenge 1006
\newsymbol\blacklozenge 1007
\newsymbol\circlearrowright 1308
\newsymbol\circlearrowleft 1309
\undefine\rightleftharpoons
\newsymbol\rightleftharpoons 130A
\newsymbol\leftrightharpoons 130B
\newsymbol\boxminus 120C
\newsymbol\Vdash 130D
\newsymbol\Vvdash 130E
\newsymbol\vDash 130F
\newsymbol\twoheadrightarrow 1310
\newsymbol\twoheadleftarrow 1311
\newsymbol\leftleftarrows 1312
\newsymbol\rightrightarrows 1313
\newsymbol\upuparrows 1314
\newsymbol\downdownarrows 1315
\newsymbol\upharpoonright 1316
 
\newsymbol\downharpoonright 1317
\newsymbol\upharpoonleft 1318
\newsymbol\downharpoonleft 1319
\newsymbol\rightarrowtail 131A
\newsymbol\leftarrowtail 131B
\newsymbol\leftrightarrows 131C
\newsymbol\rightleftarrows 131D
\newsymbol\Lsh 131E
\newsymbol\Rsh 131F
\newsymbol\rightsquigarrow 1320
\newsymbol\leftrightsquigarrow 1321
\newsymbol\looparrowleft 1322
\newsymbol\looparrowright 1323
\newsymbol\circeq 1324
\newsymbol\succsim 1325
\newsymbol\gtrsim 1326
\newsymbol\gtrapprox 1327
\newsymbol\multimap 1328
\newsymbol\therefore 1329
\newsymbol\because 132A
\newsymbol\doteqdot 132B

\newsymbol\triangleq 132C
\newsymbol\precsim 132D
\newsymbol\lesssim 132E
\newsymbol\lessapprox 132F
\newsymbol\eqslantless 1330
\newsymbol\eqslantgtr 1331
\newsymbol\curlyeqprec 1332
\newsymbol\curlyeqsucc 1333
\newsymbol\preccurlyeq 1334
\newsymbol\leqq 1335
\newsymbol\leqslant 1336
\newsymbol\lessgtr 1337
\newsymbol\backprime 1038
\newsymbol\risingdotseq 133A
\newsymbol\fallingdotseq 133B
\newsymbol\succcurlyeq 133C
\newsymbol\geqq 133D
\newsymbol\geqslant 133E
\newsymbol\gtrless 133F
\newsymbol\sqsubset 1340
\newsymbol\sqsupset 1341
\newsymbol\vartriangleright 1342
\newsymbol\vartriangleleft 1343
\newsymbol\trianglerighteq 1344
\newsymbol\trianglelefteq 1345
\newsymbol\bigstar 1046
\newsymbol\between 1347
\newsymbol\blacktriangledown 1048
\newsymbol\blacktriangleright 1349
\newsymbol\blacktriangleleft 134A
\newsymbol\vartriangle 134D
\newsymbol\blacktriangle 104E
\newsymbol\triangledown 104F
\newsymbol\eqcirc 1350
\newsymbol\lesseqgtr 1351
\newsymbol\gtreqless 1352
\newsymbol\lesseqqgtr 1353
\newsymbol\gtreqqless 1354
\newsymbol\Rrightarrow 1356
\newsymbol\Lleftarrow 1357
\newsymbol\veebar 1259
\newsymbol\barwedge 125A
\newsymbol\doublebarwedge 125B
\undefine\angle
\newsymbol\angle 105C
\newsymbol\measuredangle 105D
\newsymbol\sphericalangle 105E
\newsymbol\varpropto 135F
\newsymbol\smallsmile 1360
\newsymbol\smallfrown 1361
\newsymbol\Subset 1362
\newsymbol\Supset 1363
\newsymbol\Cup 1264

\newsymbol\Cap 1265
 
\newsymbol\curlywedge 1266
\newsymbol\curlyvee 1267
\newsymbol\leftthreetimes 1268
\newsymbol\rightthreetimes 1269
\newsymbol\subseteqq 136A
\newsymbol\supseteqq 136B
\newsymbol\bumpeq 136C
\newsymbol\Bumpeq 136D
\newsymbol\lll 136E

\newsymbol\ggg 136F

\newsymbol\circledS 1073
\newsymbol\pitchfork 1374
\newsymbol\dotplus 1275
\newsymbol\backsim 1376
\newsymbol\backsimeq 1377
\newsymbol\complement 107B
\newsymbol\intercal 127C
\newsymbol\circledcirc 127D
\newsymbol\circledast 127E
\newsymbol\circleddash 127F
\newsymbol\lvertneqq 2300
\newsymbol\gvertneqq 2301
\newsymbol\nleq 2302
\newsymbol\ngeq 2303
\newsymbol\nless 2304
\newsymbol\ngtr 2305
\newsymbol\nprec 2306
\newsymbol\nsucc 2307
\newsymbol\lneqq 2308
\newsymbol\gneqq 2309
\newsymbol\nleqslant 230A
\newsymbol\ngeqslant 230B
\newsymbol\lneq 230C
\newsymbol\gneq 230D
\newsymbol\npreceq 230E
\newsymbol\nsucceq 230F
\newsymbol\precnsim 2310
\newsymbol\succnsim 2311
\newsymbol\lnsim 2312
\newsymbol\gnsim 2313
\newsymbol\nleqq 2314
\newsymbol\ngeqq 2315
\newsymbol\precneqq 2316
\newsymbol\succneqq 2317
\newsymbol\precnapprox 2318
\newsymbol\succnapprox 2319
\newsymbol\lnapprox 231A
\newsymbol\gnapprox 231B
\newsymbol\nsim 231C
\newsymbol\ncong 231D
\newsymbol\diagup 201E
\newsymbol\diagdown 201F
\newsymbol\varsubsetneq 2320
\newsymbol\varsupsetneq 2321
\newsymbol\nsubseteqq 2322
\newsymbol\nsupseteqq 2323
\newsymbol\subsetneqq 2324
\newsymbol\supsetneqq 2325
\newsymbol\varsubsetneqq 2326
\newsymbol\varsupsetneqq 2327
\newsymbol\subsetneq 2328
\newsymbol\supsetneq 2329
\newsymbol\nsubseteq 232A
\newsymbol\nsupseteq 232B
\newsymbol\nparallel 232C
\newsymbol\nmid 232D
\newsymbol\nshortmid 232E
\newsymbol\nshortparallel 232F
\newsymbol\nvdash 2330
\newsymbol\nVdash 2331
\newsymbol\nvDash 2332
\newsymbol\nVDash 2333
\newsymbol\ntrianglerighteq 2334
\newsymbol\ntrianglelefteq 2335
\newsymbol\ntriangleleft 2336
\newsymbol\ntriangleright 2337
\newsymbol\nleftarrow 2338
\newsymbol\nrightarrow 2339
\newsymbol\nLeftarrow 233A
\newsymbol\nRightarrow 233B
\newsymbol\nLeftrightarrow 233C
\newsymbol\nleftrightarrow 233D
\newsymbol\divideontimes 223E
\newsymbol\varnothing 203F
\newsymbol\nexists 2040
\newsymbol\Finv 2060
\newsymbol\Game 2061
\newsymbol\mho 2066
\newsymbol\eth 2067
\newsymbol\eqsim 2368
\newsymbol\beth 2069
\newsymbol\gimel 206A
\newsymbol\daleth 206B
\newsymbol\lessdot 236C
\newsymbol\gtrdot 236D
\newsymbol\ltimes 226E
\newsymbol\rtimes 226F
\newsymbol\shortmid 2370
\newsymbol\shortparallel 2371
\newsymbol\smallsetminus 2272
\newsymbol\thicksim 2373
\newsymbol\thickapprox 2374
\newsymbol\approxeq 2375
\newsymbol\succapprox 2376
\newsymbol\precapprox 2377
\newsymbol\curvearrowleft 2378
\newsymbol\curvearrowright 2379
\newsymbol\digamma 207A
\newsymbol\varkappa 207B
\newsymbol\Bbbk 207C
\newsymbol\hslash 207D
\undefine\hbar
\newsymbol\hbar 207E
\newsymbol\backepsilon 237F

%
\newcount\marknumber	\marknumber=1
\newcount\countdp \newcount\countwd \newcount\countht
%
%
\ifx\pdfoutput\undefined
\def\rgboo#1{}
\def\postscript#1{\special{" #1}}		
\postscript{
	/bd {bind def} bind def
	/fsd {findfont exch scalefont def} bd
	/sms {setfont moveto show} bd
	/ms {moveto show} bd
	/pdfmark where		
	{pop} {userdict /pdfmark /cleartomark load put} ifelse
	[ /PageMode /UseOutlines		
	/DOCVIEW pdfmark}
\def\bookmark#1#2{\postscript{		
	[ /Dest /MyDest\the\marknumber /View [ /XYZ null null null ] /DEST pdfmark
	[ /Title (#2) /Count #1 /Dest /MyDest\the\marknumber /OUT pdfmark}%
	\advance\marknumber by1}
\def\pdfclink#1#2#3{%
	\hskip-.25em\setbox0=\hbox{#2}%
		\countdp=\dp0 \countwd=\wd0 \countht=\ht0%
		\divide\countdp by65536 \divide\countwd by65536%
			\divide\countht by65536%
		\advance\countdp by1 \advance\countwd by1%
			\advance\countht by1%
		\def\linkdp{\the\countdp} \def\linkwd{\the\countwd}%
			\def\linkht{\the\countht}%
	\postscript{
		[ /Rect [ -1.5 -\linkdp.0 0\linkwd.0 0\linkht.5 ]
		/Border [ 0 0 0 ]
		/Action << /Subtype /URI /URI (#3) >>
		/Subtype /Link
		/ANN pdfmark}{\rgb{#1}{#2}}}
%
%
\else
\def\rgboo#1{\pdfliteral{#1 rg #1 RG}}
\pdfcatalog{/PageMode /UseOutlines}		
\def\bookmark#1#2{
	\pdfdest num \marknumber xyz
	\pdfoutline goto num \marknumber count #1 {#2}
	\advance\marknumber by1}
\def\pdfklink#1#2{%
	\noindent\pdfstartlink user
		{/Subtype /Link
		/Border [ 0 0 0 ]
		/A << /S /URI /URI (#2) >>}{\rgb{1 0 0}{#1}}%
	\pdfendlink}
\fi

\def\rgbo#1#2{\rgboo{#1}#2\rgboo{0 0 0}}
\def\rgb#1#2{\mark{#1}\rgbo{#1}{#2}\mark{0 0 0}}
\def\pdflink#1{\pdfklink{#1}{#1}}
%
%
\newcount\seccount  
\newcount\subcount  
\newcount\ssscount  
\newcount\clmcount  
\newcount\equcount  
\newcount\refcount  
\newcount\demcount  
\newcount\execount  
\newcount\procount  
\seccount=0
\equcount=0
\clmcount=0
\subcount=0
\refcount=0
\demcount=0
\execount=0
\procount=0
%
\def\proof{\medskip\noindent{\bf Proof.\ }}
\def\proofof(#1){\medskip\noindent{\bf Proof of \csname c#1\endcsname.\ }}
\def\qed{\hfill{\sevenbf QED}\par\medskip}
\def\references{\bigskip\noindent\hbox{\bf References}\medskip
                \ifx\pdflink\undefined\else\bookmark{0}{References}\fi}
\def\addref#1{\global\advance\refcount by 1
              \expandafter\xdef\csname r#1\endcsname{\number\refcount}}

\def\nextremark #1\par{\item{$\circ$} #1}
\def\firstremark #1\par{\bigskip\noindent{\bf Remarks.}
     \smallskip\nextremark #1\par}
\def\abstract#1\par{{\baselineskip=10pt
    \eightpoint\narrower\noindent{\eightbf Abstract.} #1\par}}
%
\def\equtag#1{\global\advance\equcount by 1
              \expandafter\xdef\csname e#1\endcsname{(\number\seccount.\number\equcount)}}
\def\equation(#1){\equtag{#1}\eqno\csname e#1\endcsname}
\def\equ(#1){\hskip-0.03em\csname e#1\endcsname}
%
\def\clmtag#1#2{\global\advance\clmcount by 1
                \expandafter\xdef\csname cn#2\endcsname{\number\seccount.\number\clmcount}
                \expandafter\xdef\csname c#2\endcsname{#1~\number\seccount.\number\clmcount}}
\def\claim #1(#2) #3\par{\clmtag{#1}{#2}
    \vskip.1in\medbreak\noindent
    {\bf \csname c#2\endcsname .\ }{\sl #3}\par
    \ifdim\lastskip<\medskipamount
    \removelastskip\penalty55\medskip\fi}
\def\clm(#1){\csname c#1\endcsname}
\def\clmno(#1){\csname cn#1\endcsname}
%
\def\sectag#1{\global\advance\seccount by 1
              \expandafter\xdef\csname sectionname\endcsname{\number\seccount. #1}
              \equcount=0 \clmcount=0 \subcount=0 \execount=0 \procount=0}
\def\section#1\par{\vskip0pt plus.1\vsize\penalty-40
    \vskip0pt plus -.1\vsize\bigskip\bigskip
    \sectag{#1}
    \message{\sectionname}\leftline{\twelvebf\sectionname}
    \ifx\pdflink\undefined
    \else
      \bookmark{0}{\sectionname}
    \fi
    \nobreak\smallskip\noindent}
%
\def\subtag#1{\global\advance\subcount by 1
              \expandafter\xdef\csname subsectionname\endcsname{\number\seccount.\number\subcount. #1}
              \ssscount=0}
\def\subsection#1\par{\vskip0pt plus.05\vsize\penalty-20
    \vskip0pt plus -.05\vsize\medskip\medskip
    \subtag{#1}
    \message{\subsectionname}\leftline{\tenbf\subsectionname}
    \ifx\pdflink\undefined
    \else
      \bookmark{0}{... \subsectionname}  
    \fi
    \nobreak\smallskip\noindent}
%
\def\ssstag#1{\global\advance\ssscount by 1
              \expandafter\xdef\csname subsubsectionname\endcsname{\number\seccount.\number\subcount.\number\ssscount. #1}}
\def\subsubsection#1\par{\vskip0pt plus.02\vsize\penalty-10
    \vskip0pt plus -.02\vsize\smallskip\smallskip
    \ssstag{#1}
    \message{\subsubsectionname}\leftline{\ninebf\subsubsectionname}
    \nobreak\smallskip\noindent
    \ifx\pdflink\undefined
    \else
      \bookmark{0}{...... \subsubsectionname}  
    \fi}
%
\def\demtag#1#2{\global\advance\demcount by 1
                \expandafter\xdef\csname de#2\endcsname{#1~\number\demcount}}
\def\demo #1(#2) #3\par{
  \demtag{#1}{#2}
  \vskip.1in\medbreak\noindent
  {\bf #1 \number\demcount.\enspace}
  {\rm #3}\par
  \ifdim\lastskip<\medskipamount
  \removelastskip\penalty55\medskip\fi}
\def\dem(#1){\csname de#1\endcsname}
%
\def\exetag#1{\global\advance\execount by 1
              \expandafter\xdef\csname ex#1\endcsname{Exercise~\number\seccount.\number\execount}}
\def\exercise(#1) #2\par{
  \exetag{#1}
  \vskip.1in\medbreak\noindent
  {\bf Exercise \number\execount.}
  {\rm #2}\par
  \ifdim\lastskip<\medskipamount
  \removelastskip\penalty55\medskip\fi}
\def\exe(#1){\csname ex#1\endcsname}
%
\def\protag#1{\global\advance\procount by 1
              \expandafter\xdef\csname pr#1\endcsname{\number\seccount.\number\procount}}
\def\problem(#1) #2\par{
  \ifnum\procount=0
    \parskip=6pt
    \vbox{\bigskip\centerline{\bf Problems \number\seccount}\nobreak\medskip}
  \fi
  \protag{#1}
  \item{\number\procount.} #2}
\def\pro(#1){Problem \csname pr#1\endcsname}
%
%
%
\def\rightheadline{\hfil}
\def\leftheadline{\sevenrm\hfil HANS KOCH\hfil}
\headline={\ifnum\pageno=\firstpage\hfil\else
\ifodd\pageno{{\fiverm\rightheadline}\number\pageno}
\else{\number\pageno\fiverm\leftheadline}\fi\fi}
\footline={\ifnum\pageno=\firstpage\hss\tenrm\folio\hss\else\hss\fi}

\let\cl=\centerline

\let\eps=\varepsilon
\let\sss=\scriptscriptstyle

\def\AA{{\cal A}}
\def\BB{{\cal B}}
\def\CC{{\cal C}}

\def\FF{{\cal F}}
\def\GG{{\cal G}}
\def\HH{{\cal H}}

\def\LL{{\cal L}}

\def\OO{{\cal O}}
\def\PP{{\cal P}}

\def\RR{{\cal R}}
\def\SS{{\cal S}}

\def\WW{{\cal W}}

\def\ZZ{{\cal Z}}
\def\ssA{{\sss A}}
\def\ssB{{\sss B}}

\def\ssE{{\sss E}}

\def\ssL{{\sss L}}

\def\ssN{{\sss N}}

\def\id{{\rm I}}

\def\tr{\mathop{\rm tr}\nolimits}
\def\det{\mathop{\rm det}\nolimits}

\def\Im{\mathop{\rm Im}\nolimits}
%
\newfam\dsfam
\def\mathds #1{{\fam\dsfam\tends #1}}

\font\tends=dsrom10
\font\eightds=dsrom8
\textfont\dsfam=\tends
\scriptfont\dsfam=\eightds
%

\def\integer{{\mathds Z}}
\def\rational{{\mathds Q}}
\def\real{{\mathds R}}
\def\complex{{\mathds C}}

\def\torus{{\Bbb T}}

\def\bdot{\hbox{\bf .}}
\def\bcomma{\hbox{\bf ,}}
\def\defeq{\mathrel{\mathop=^{\sss\rm def}}}
\def\half{{1\over 2}}

\def\thalf{{\textstyle\half}}

\def\twovec#1#2{\left[\matrix{#1\cr#2\cr}\right]}

\def\twomat#1#2#3#4{\left[\matrix{#1&#2\cr#3&#4\cr}\right]}

%

%

%

%

\font\tenamsb=msbm10 \font\sevenamsb=msbm7 \font\fiveamsb=msbm5
\newfam\bbfam
\textfont\bbfam=\tenamsb
\scriptfont\bbfam=\sevenamsb
\scriptscriptfont\bbfam=\fiveamsb

\def\ttA{{\tt A}}
\def\ttB{{\tt B}}

\def\firstAM{{\ninerm AM}}
\def\AM{{\ninerm AM~}}
\def\firstRG{{\ninerm RG}}
\def\RG{{\ninerm RG~}}

\def\onedim{\hbox{\rm 1d}}
\def\twodim{\hbox{\rm 2d}}
\def\rad{r}
\def\radB{{\rad_{_{\hskip-1pt B}}}}
\def\radA{{\rad_{_{\hskip-1pt A}}}}
\def\bua{{\hbox{\magnineeufm a}}}
\def\bub{{\hbox{\magnineeufm b}}}
\def\buk{{\hbox{\magnineeufm k}}}
\def\buN{{\hbox{\teneufm N}}}
\def\buR{{\hbox{\teneufm R}}}

\font\tenib=cmmib10
\def\bma{{\hbox{\tenib a}}}

%

%
\def\circle{{\Bbb S}}

\def\Rot{\mathop{\rm Rot}\nolimits}
\def\rot{\mathop{\rm rot}\nolimits}
\def\mod{\mathop{\rm mod}\nolimits}
\def\gcd{\mathop{\rm gcd}\nolimits}
\def\idmat{{\bf 1}}
\def\rmGL{{\rm GL}}
\def\rmSL{{\rm SL}}
\def\ssA{{\sss A}}
\def\ssB{{\sss B}}

\def\symm{{\sss\circ}}

\def\tfrac#1#2{{\textstyle{#1\over #2}}}
\def\bdot{\hbox{\bf .}}
\def\hdots{\line{\leaders\hbox to 0.5em{\hss .\hss}\hfil}}

\def\sfrac#1#2{\hbox{\raise2.2pt\hbox{$\scriptstyle#1$}\hskip-1.2pt
   {$\scriptstyle/$}\hskip-0.9pt\lower2.2pt\hbox{$\scriptstyle#2$}\hskip1.0pt}}
\def\shalf{\sfrac{1}{2}}
\def\squarter{\sfrac{1}{4}}

\def\stwomat#1#2#3#4{{\eightpoint\left[\matrix{#1&#2\cr#3&#4\cr}\right]}}

\def\today{\ifcase\month\or
January\or February\or March\or April\or May\or June\or
July\or August\or September\or October\or November\or December\fi
\space\number\day, \number\year}
\addref{Harp}
\addref{Sos}
\addref{Sur}
\addref{Swi}
\addref{Kato}
\addref{Hof}
\addref{TKNdN}
\addref{JoMo}
\addref{Herman}
\addref{DelSou}
\addref{AvSi}
\addref{Johnson}
\addref{WieZa}
\addref{BeGap}
\addref{Lastiii}
\addref{FaKa}
\addref{HKW}
\addref{Jito}
\addref{BouJitoii}
\addref{Puigi}
\addref{AvKri}
\addref{Dama}
\addref{GoSch}
\addref{JiMarx}
\addref{ATZ}
\addref{KochAM}
\addref{KochTrig}
\addref{KochUniv}
\def\leftheadline{\sixrm\hfil H.~Koch\hfil\today}
\def\rightheadline{\sevenrm\hfil Supercritical skew-product maps\hfil}
%
\cl{{\twelvebf Universal supercritical behavior for some skew-product maps}}
\bigskip

\cl{
Hans Koch
\footnote{$\!^1$}
{\eightpoint\hskip-2.7em
Department of Mathematics, The University of Texas at Austin,
Austin, TX 78712.}}

\bigskip
\abstract
We consider skew-product maps over circle rotations $x\mapsto x+\alpha\;(\mod\,1)$
with factors that take values in $\rmSL(2,{\scriptstyle\real})$.
In numerical experiments with $\alpha$ the inverse golden mean,
Fibonacci iterates of almost Mathieu maps
with rotation number $\scriptstyle 1/4$ and positive Lyapunov exponent
exhibit asymptotic scaling behavior.
We prove the existence of such asymptotic scaling
for ``periodic'' rotation numbers and for large Lyapunov exponent.
The phenomenon is universal,
in the sense that it holds for open sets of maps,
with the scaling limit being independent of the maps.
The set of maps with a given periodic rotation number
is a real analytic codimension $1$ manifold in a suitable space of maps.

\section Introduction and main results

We consider the asymptotic behavior of skew products
$$
\eqalign{
A^{\ast q}_\symm(x)
&=A^{\ast q}\bigl(x-\tfrac{q}{2}\alpha\bigr)\,,\cr
A^{\ast q}(x)
&\,\defeq\,
A(x+(q-1)\alpha)\cdots A(x+2\alpha)A(x+\alpha)A(x)\,,\cr}
\equation(Aastqx)
$$
as $q\to\infty$ along certain subsequences,
where $\alpha$ is an irrational number
and $A$ is a real analytic function
from the circle $\torus=\real/\integer$ to the group $\rmSL(2,\real)$.
Products for negative $q$ are defined
by replacing the factors $A$ in the above equation by $A(\,\bdot-\alpha)^{-1}$.

The typical growth of such products is described
by the Lyapunov exponent.
The Lyapunov exponent of a pair $G=(\alpha,A)$
with $A:\torus\to\rmSL(2,\real)$ continuous is defined as
$$
L(G)=\lim_{q\to\infty}{1\over q}\log\bigl\|A^{\ast q}(x)\bigr\|\,.
\equation(LyapDef)
$$
Assuming that $\alpha$ is irrational,
this limit exists by ergodicity and is a.e.~constant in $x$.

\smallskip
The work presented in this paper was motivated in part
by the observation in [\rKochUniv] of systematic
slower-than-typical growth of certain products, as will be described below.
We restrict to factors that are reversible in the following sense.

\claim Definition(reversibleFac)
The symmetric factor $A_\symm$ associated with a pair $G=(\alpha,A)$
is defined by setting $A_\symm(x)=A\bigl(x-\sfrac{\alpha}{2}\bigr)$ for all $x$.
We say that $G$ is reversible if $S^{-1}A_\symm(x)S=A_\symm(-x)^\dagger$ for all $x$,
where
$$
S=\stwomat{0}{1}{1}{0}\,,\qquad
M^\dagger=\stwomat{d}{-b}{-c}{a}\quad{\rm if}\quad
M=\stwomat{a}{b}{c}{d}\,.
\equation(reversibleFac)
$$
Notice that $M^\dagger M=\det(M)\idmat$.
We will refer to $M^\dagger$ as the quasi-inverse of $M$.

To be more precise about the type of products being considered,
let $\alpha={\sqrt{5}-1\over 2}$ be the inverse golden mean
and $p_n$ the $n$-th Fibonacci number.
Then the following holds.

\claim Theorem(subseqLimit)
Let $A:\torus\to\rmSL(2,\real)$ be real analytic and reversible.
Assume that $G=(\alpha,A)$ has a positive Lyapunov exponent.
Then there exist real numbers $M_1,M_2,M_3,\ldots$
such that every subsequence of
$$
n\mapsto M_n\,A^{\ast p_n}_\symm\bigl(\alpha^n\bdot\bigr)
\equation(SymmFacSeq)
$$
has a subsequence that converges to a nonzero function $z\mapsto b_\diamond(z)\ttB$,
where $\ttB$ is a constant $2\times 2$ matrix of rank one.
The scalar factor $b_\diamond$ is an even entire function,
and convergence is uniform on compact subsets of $\complex$.

A more restrictive version of this theorem was given in [\rKochUniv].
The main objects under investigation in [\rKochUniv] were parametrized
factors $A$ at the boundary between zero Lyapunov exponent
(subcritical or critical factors) and positive Lyapunov exponent (supercritical factors).
In this context, products of supercritical factors were computed numerically,
and limit functions $b_\diamond$ were found with surprisingly regular zeros.
Our goal here is to describe this type of limits.

Our main focus is on Schr\"odinger factors
$$
A(x)=\twomat{\lambda v(t+x)-E}{-1}{1}{0}\,,
\equation(StdAM)
$$
and perturbations of such factors.
Here, $\lambda$, $t$, and $E$ are real parameters,
and $v$ is a $1$-periodic continuous function on $\real$.
The choice $v(x)=-2\cos(2\pi x)$ defines
the almost Mathieu (\firstAM) family of factors [\rLastiii,\rDama,\rJiMarx].
Schr\"odinger factors arise naturally in the study of Schr\"odinger operators
on $\ell^2(\integer)$ of the form
\footnote{$\!^2$}
{\eightpoint\hskip-22pt
We adopt here the ``physical'' choice of signs,
which makes the kinetic term positive.}
$$
(\HH_\lambda^t u)_q=\lambda v(t+q\alpha)u_q-u_{q+1}-u_{q-1}\,,\qquad q\in\integer\,.
\equation(origAMHam)
$$
To be more specific, $A^{\ast q}(0)$ is the transfer matrix that maps
$\bigl[{u_0\atop u_{-1}}\bigr]$ to $\bigl[{u_q\atop u_{q-1}}\bigr]$,
for a formal solution $u$ of the equation $\HH_\lambda^t u=Eu$.
Many properties of $\HH_\lambda^t$ are reflected
by properties of the corresponding products \equ(Aastqx), and vice versa.
The Lyapunov exponent of the \AM factors with $\lambda>0$
is known to be $L(G)=\max\{0,\log\lambda\}$ for all energies
in the spectrum [\rBouJitoii].

The family of \AM operators $t\mapsto\HH_\lambda^t$
describe the motion of an electron on $\integer^2$
under the influence of a magnetic flux $2\pi\alpha$ per unit cell,
after restricting to wave functions $\phi(q,p)=e^{-2\pi ip t}u_q$.
The full Hamiltonian for this system
is known as the Hofstadter Hamiltonian [\rHarp,\rHof].

\smallskip
Another quantity associated with $\HH_\lambda^t$ that is relevant in the present context
is the integrated density of states $E\mapsto\buk(E)$.
Notice that $\HH_\lambda^t$ is self-adjoint and bounded.
Denote by $P_\ssL$ the canonical inclusion map
from $\ell^2\bigl(\integer\cap[-L,L]\bigr)$ into $\ell^2(\integer)$.
Given $E\in\real$,
the integrated density of states $\buk_\ssL(\HH_\lambda^t,E)$ is defined as the fraction
of eigenvalues of $P_\ssL^\ast\HH_\lambda^t P_\ssL$ that belong to $(-\infty,E]$.
The integrated density of states for $\HH_\lambda^t$ can be obtained
as the limit $\buk(E)=\lim_{L\to\infty}\buk_\ssL(\HH_\lambda^t,E)$.
Clearly, $\buk$ is an increasing function,
taking the value $\buk(E)=0$ for $E$ below the spectrum of $\HH_\lambda^t$
and the value $\buk(E)=1$ for $E$ above the spectrum.
And it is constant on any spectral gaps
(connected components of the resolvent set).
For the \AM operator, it is known that there are no isolated eigenvalues [\rAvKri],
so $\buk$ is continuous.
For $\lambda>1$, the spectrum is pure point,
with exponentially decreasing eigenfunctions [\rJito].
\big(This holds for arbitrary Diophantine $\alpha$,
and for ``nonresonant'' values of $t$
that include our choice $\sfrac{\alpha}{2}$ below.\big)

Among the intriguing features of the \AM operator
is the following resonance phenomenon.
$\HH^t_\lambda$ has an infinite number of spectral gaps.
Each gap can be labeled canonically by an integer $k$,
known as the Hall conductance,
and $\buk(E)\equiv k\alpha\;(\mod 1)$ holds for all energies $E$
in the spectral gap labeled by $k$.
For details we refer to [\rTKNdN,\rJoMo,\rAvSi,\rBeGap,\rGoSch].

\smallskip
A related quantity for Schr\"odinger pairs $G=(\alpha,A)$ is the
rotation number $\rot(G)$.
In cases where $A$ is of the form \equ(StdAM),
with a potential $v$ that is continuous,
we can define the rotation number by the equation
$$
\rot(G)=\lim_{q\to\infty}{\Rot_\ssN(G)\over N}\,,
\qquad\qquad\hbox{\rm (for Schr\"odinger $G$)}
\equation(rotviaSigns)
$$
where $\Rot_\ssN(G)$ is half the number of sign changes
on $\{1,2,\ldots,N\}$
of a nonzero solution $u$ of the equation $\HH_\lambda^t u=Eu$,
adding $\shalf$ if $u(N)=0$.
Notice that $0\le\rot(G)\le\shalf$.
The rotation number for more general pairs $G$,
and its properties, will be described in Section 3.
Assuming that $\alpha$ is irrational,
the value of $\rot(G)$ is independent of the solution $u$.
For an \AM pair $G_\ssE$ with spectral energy $E$,
it can be shown [\rAvSi] that the rotation number
is related to the integrated density of states via $\buk(E)=2\rot(G_\ssE)$.
Roughly speaking, this relation expresses the fact
that the wave number (inverse wave length)
of eigenfunctions increases with the energy $E$.

\demo Convention(FixPhase)
In what follows, we fix the phase $t$ in the definition
of the Schr\"odinger factor \equ(StdAM)
and the Schr\"odinger operator \equ(origAMHam)
to the value $t=\sfrac{\alpha}{2}$, unless specified otherwise.
This makes the symmetric \AM factors $A_\symm$ reversible
in the sense of \clm(reversibleFac).

The above-mentioned gap-labeling property
associates spectral gaps with rotation numbers
that belong to $\half\integer[\alpha]$.
We will now describe a different phenomenon
that is associated with some other rotation numbers in $\rational[\alpha]$.

\claim Definition(ResonantDef)
A number $\rho\in[0,\shalf]$ will be called positive periodic
(as a rotation number) if
$\rho={w\over v}+{u\over v}\alpha$ for some integers $u,v,w$
with $v>2$ and $\gcd(u,v,w)=1$,
and if $\left|{u\over v}-{w\over v}\alpha\right|\le\half$.

We will call $\rho\in[-\shalf,0]$ negative periodic if $\rho+\shalf$ is positive periodic.
Notice that, for any given integer $v>2$, the (positive or negative)
periodic rotation numbers in ${1\over v}\integer[\alpha]$
constitute a discrete subset of $\real$.

In what follows, $\alpha$ denotes the inverse
golden mean, unless specified otherwise.
Let $p_n$ be the $n$-th Fibonacci number, and define $q_n=p_{n+1}$.
Our main result in this paper is the following.

\claim Theorem(supercritAMLimits)
Let $\rho$ be positive periodic.
Then there exists a positive integer $n$
and two nonzero even entire functions $b_\diamond$ and $a_\diamond$,
such that the following holds.
Consider the family of \AM maps $G_{\lambda,E}$
parametrized by the coupling constant $\lambda>0$ and the energy $E$.
Then there exists a real analytic function $\epsilon$,
defined on an open neighborhood of zero,
such that if $\delta=\lambda^{-1}$ belongs to the domain of $\epsilon$,
and if $E=\lambda\epsilon(\delta)$,
then $G$ has rotation number $\rot(G)=\rho$.
In this case, there exist sequences of positive real numbers
$k\mapsto M_k$ and $k\mapsto W_k$, such that
$$
\eqalign{
\lim_{k\to\infty}M_k\,A^{\ast p_{3nk}}_\symm\bigl(\alpha^{3nk}x\bigr)
&=b_\diamond(x)\ttA^\dagger\,,\cr
\lim_{k\to\infty}W_k\,A^{\ast q_{3nk}}_\symm\bigl(\alpha^{3nk}x\bigr)
&=a_\diamond(x)\ttA\,,\cr}
\equation(rhoSymmFacSeq)
$$
for some constant $2\times 2$ matrix $\ttA$ of rank one.
Convergence in \equ(rhoSymmFacSeq) is uniform on compact
subsets of $\complex$.
An analogous result holds for any two-parameter family
that is sufficiently close
to the \AM family (in a sense described later).

Both $b_\diamond$ and $a_\diamond$ have infinitely many zeros,
all simple and on the real axis.
The zeros can be determined rather explicitly,
and this is used to construct the functions $b_\diamond$ and $a_\diamond$.
We note that the function $\epsilon=\epsilon(\delta)$
and the matrix $\ttA$ in the above theorem can depend on the family.
But the scaling limits $b_\diamond$ and $a_\diamond$ are universal.

We expect that an analogous result holds for all $\lambda>1$.
Numerically, this is observed at energy $E=0$,
which corresponds to $\rho=\squarter$.
For this particular value of $\rho$, a proof seems feasible,
since the emergence of the limiting zeros
is quite transparent in this case [\rKochUniv].
We also expect that \clm(supercritAMLimits)
generalizes to arbitrary quadratic irrationals $\alpha$.

\demo Remark(PosLambda)
For non-Schr\"odinger maps $G=(\alpha,A)$,
the rotation number $\rot(G)$ in \clm(supercritAMLimits)
has to be replaced by the rotation number $\varrho(G)$ defined in Section 3.
Using $\varrho(G)$ in place of $\rot(G)$,
an analogous theorem holds for negative periodic rotation numbers as well.
This corresponds to replacing $A$ by $-A$.
For a Schr\"odinger factor $A$,
we can get $-A$ via a conjugacy by $\bigl[{1~\phantom{-}0\atop 0~-1}\bigr]$
and then replacing $v,E$ by their negatives.
The conjugacy is equivalent to replacing $u$ by $q\mapsto(-1)^qu_q$
in the equation \equ(origAMHam).

\clm(supercritAMLimits) is best understood in a dynamics context.
The pairs $G=(\alpha,A)$ with $A(x)\in\rmSL(2,\real)$
are naturally associated with skew-product maps
$$
G(x,y)=\bigl(x+\alpha,A(x)y\bigr)\,,\qquad x\in X\,,\quad y\in Y\,,
\equation(MapGDef)
$$
on $X\times Y$, where $X=\torus$ and $Y=\real^2$.
The $q$-th iterate of $G$ is $G^q=\bigl(q\alpha,A^{\ast q})$,
with $A^{\ast q}$ as defined in \equ(Aastqx).
So \clm(supercritAMLimits) can be viewed as describing
the asymptotic behavior of long orbits for skew-product maps.
Due to the scaling involved in \equ(rhoSymmFacSeq),
we will also need to consider the case $X=\real$.
The condition $\det(A)=1$ will be weakened as well,
and $Y=\real$ will be used for limit cases.

Our proof of \clm(supercritAMLimits) involves the use
of a renormalization transformation $\RR$ that acts on
pairs $P=(F,G)$ of skew-product maps $F=(1,B)$ and $G=(\alpha,A)$.
We consider families of pairs $P_{\delta,\epsilon}$
that admit a scaling limit as $\delta\to 0$.
For Schr\"odinger factors \equ(StdAM), we use $\delta=\lambda^{-1}$,
multiply $A$ by $\delta$, and set $E=\lambda\epsilon$.
If $\varrho(G_{\delta,\epsilon})=\rho$ with $\rho$ positive periodic,
then the asymptotic behavior \equ(rhoSymmFacSeq)
is governed by a fixed point $P_\diamond$ of
(a modified version of) the transformation $\RR^{3n}$.
The local stable manifold of this transformation at $P_\diamond$
is of codimension $1$ and agrees with the level set $\varrho(G)=\rho$.
So the real analytic curve $\delta\mapsto\epsilon(\delta)$
described in \clm(supercritAMLimits)
characterizes the intersection of the given family with this manifold.

Given the relation $2\rot(G_\ssE)=\buk(E)$ for the \AM family,
the above suggests that the energies $E=\lambda\epsilon\bigl(\lambda^{-1}\bigr)$
are eigenvalues of the operator $\HH^{\alpha/2}_\lambda$.
Analytic curves of eigenvalues are a standard feature
in real analytic families of compact self-adjoint operators,
even in the presence of level crossings [\rKato].
This applies e.g.~to the operators $P_\ssL^\ast\HH^t_\lambda P_\ssL$
mentioned earlier.
So it seems possible that $E=\lambda\epsilon\bigl(\lambda^{-1}\bigr)$
is a limit as $L\to\infty$ of such eigenvalue curves.

{}From a dynamics perspective,
the zeros at $\pm\rho$ of the limit function $a_\diamond$
seem associated with factors $A_\symm(x)$ in large products
that map the expanding direction
(from the product to the right of this factor)
to the contracting direction of the product to the left of this factor.
This suggest that $x$ is near $\pm\rho$, which is
where $A_\symm(x)$ can cover a wide range of rotation angles.
The location of such factors
may be associated with peaks of an eigenvector.

We note that the Hofstadter model has
a large number of symmetries [\rWieZa,\rFaKa,\rHKW]
and other important arithmetic features.
It seems fair to say that their interplay
with analysis is rather poorly understood.
This is our motivation for focusing on the inverse golden mean.
Periodicity may be regarded as exceptional,
but periodic orbits play an important role in dynamical systems,
since they heavily influence the motion nearby.
The \AM family can be regarded as exceptional as well,
but, as \clm(supercritAMLimits) shows,
some of its ``rigid'' behavior is shared by nearby families.

\smallskip
The remaining part of this paper is organized as follows.
In Section 2 we introduce the renormalization transformation $\RR$
and prove a result that implies \clm(subseqLimit).
Section 3 is devoted to limits of Schr\"odinger factors
and the rotation number.
In Section 4 we construct the periodic orbits of $\RR$
that are associated with positive periodic rotation numbers.
A proof of \clm(supercritAMLimits) is given in Section 5.

\demo Remarks(Some)

\nextremark
To leading order, the scaling constants in \clm(supercritAMLimits)
are $M_k\sim e^{-p_{3nk}L}$ and $W_k\sim e^{-q_{3nk}L}$,
where $L$ is the Lyapunov exponent of $G$.
For the \AM factors, we suspect that a choice
$M_k=Me^{-q_{3nk}L}$ and $W_k=We^{-p_{3nk}L}$ will work.
But trying to verify this would go beyond the scope of this paper.

\nextremark
If we allowed $v=1$ and $v=2$ in \clm(ResonantDef),
then $\{0,\sfrac{\alpha}{2},\shalf\}$
would count as positive periodic as well.
The associated energies for the \AM operators
correspond to spectral gaps.
We believe that \clm(supercritAMLimits) can be proved
for these rotation numbers as well,
but they would require special treatment.

\nextremark
Numerically, critical fixed points
\big(for the appropriate power of $\RR^3$\big)
have been found for rotation numbers in
$\{0,\sfrac{1}{8},\sfrac{1}{6},\shalf-\sfrac{\alpha}{2},
\sfrac{1}{4},\sfrac{\alpha}{2},\sfrac{2}{6},\sfrac{3}{8},\shalf\}$.
Existence proofs have been given in [\rKochAM,\rKochUniv]
in the case of rotation numbers $\{0,\sfrac{1}{4},\shalf\}$.

\nextremark
Conjecture 1.1 in [\rKochUniv] on critical periods
only covers rational rotation numbers.
A reasonable amendment would be to include all (positive or negative)
periodic rotation numbers.

\nextremark
Numerically, the unstable manifold
associated with the critical fixed point of $\RR^3$
described in Theorem 2.3 of [\rKochUniv]
leads into a supercritical fixed point $P_\diamond$
of the type described here.
This concerns the case $\rho=\squarter$,
where the curve $E=\lambda\epsilon\bigl(\lambda^{-1}\bigr)$ is trivial,
namely $E=0$ by symmetry.
Still, a possible renormalization picture would be that
something similar occurs for any periodic rotation number $\rho$.

\section Renormalization

We use renormalization as a combinatorial tool
for generating factors of the type that appear in \equ(rhoSymmFacSeq).
{}From an arithmetic point of view,
the renormalization (\firstRG) transformation $\RR$ defined below
lifts the Gauss map for real numbers to pairs of maps $(F,G)$.
From an analysis point of view,
$\RR$ constitutes a dynamical system on a space of such pairs.

Consider a pair $(F,G)$ of skew-product maps
$F=(1,B)$ and $G=(\alpha,A)$ on $\real\times\real^2$,
with factors $B,A:\real\to\rmSL(2,\real)$.
If $G$ commutes with $F$, then $G$ can be viewed
as a map on a cylinder in $\real\times\real^2$
whose points are the orbits of $F$.
Assume that $\alpha<1$ is a positive irrational number,
and denote by $c$ the integer part of $\alpha^{-1}$.
Then the basic renormalized pair is defined by
$$
\RR(F,G)=\bigl(\check F,\check G\bigr)\,,\qquad
\check F=\Lambda^{-1}G\Lambda\,,\quad
\check G=\Lambda^{-1}FG^{-c}\Lambda\,,
\equation(RGDef)
$$
where $\Lambda(x,y)=(\alpha x,y)$.
By construction, the first component of $\check F$ is again $1$,
while the first component of $\check G$ is $\check\alpha=\alpha^{-1}-c$.
Notice that $\alpha\mapsto\check\alpha$ is the Gauss map
that appears in the continued fraction expansion of $\alpha$.

In what follows, we restrict to the inverse golden mean,
which satisfies $\alpha^{-1}-1=\alpha$.
So we can fix $c=1$ in \equ(RGDef).
The transformation $\RR$ extends readily to factors in $\rmGL(2,\real)$.
In this case, $G^{-1}$ is replaced by the quasi-inverse $G^\dagger$ defined below.
We note that a map $G=(\alpha,A)$ is reversible in the sense of \clm(reversibleFac)
precisely if
$$
\SS^{-1}G\SS=G^\dagger\,\defeq\,\bigl(-\alpha,A(\,\bdot-\alpha)^\dagger\bigr)\,,\qquad
\SS(x,y)=(-x,Sy)\,.
\equation(mapReversible)
$$
A pair $P=(F,G)$ is said to be reversible if both $F$ and $G$ are reversible.

The conjugacy by $\Lambda$ only renormalizes the translational
parts of our skew-products.
We will often have to renormalize the matrix part as well.
When considering periodic orbits of length $k$,
it is best to do this only once per period.
So for $\RR^k$ we use a scaling
$$
\Lambda_k=\LL_k\Lambda^k\,,\qquad
\LL_k(x,y)=(x,L_k y)\,,
\equation(LLkDef)
$$
where $L_k$ a suitable matrix in $\pm\rmSL(2,\real)$ that commutes with $S$.
This matrix can depend on the pair $P=(F,G)$.
It will be specified as needed, when the choice matters.

\smallskip
The definition \equ(RGDef) of $\RR$ involves
the basic composition operator $\CC(F,G)=(G,FG^\dagger)$.
Due to the symmetry $A\mapsto-A$ mentioned in \dem(PosLambda),
we restrict to powers of $\CC$ that are multiples of $3$.
Then it is convenient to replace $\CC^3$ by the transformation
$\CC_3$ defined by
$$
\CC_3(F,G)=\bigl(GF^\dagger G\,\bcomma\,G^\dagger FG^\dagger FG^\dagger\bigr)\,.
\equation(CCThree)
$$
The only difference between $\CC^3$ and $\CC_3$ is the order of the factors.
This is irrelevant for commuting pairs;
but for non-commuting pairs, which need to be included in our analysis,
the transformation $\CC^3$ does not in general preserve reversibility,
while $\CC_3$ does.

Given any $n\ge 1$, we define the \RG transformation $\RR_{3n}$
by setting
$$
\RR_{3n}(F,G)=\bigl(\Lambda_{3n}^{-1}\hat F\Lambda_{3n}\,\bcomma\,
\Lambda_{3n}^{-1}\hat F\Lambda_{3n}\bigr)\,,
\qquad\bigl(\hat F,\hat G\bigr)=\CC_3^n(F,G)\,.
\equation(PaliRG)
$$

Consider first the case $n=1$.
Let $\bigl(\tilde F,\tilde G\bigr)=\RR_3(F,G)$.
A straightforward computation shows that
the symmetric factor $\tilde B_\symm$ for $\tilde F$ is
$$
\tilde B_\symm(x)
=L_3^{-1}A_\symm\bigl(\alpha^3 x-\tfrac{1-\alpha}{2}\bigr)
B_\symm\bigl(\alpha^3 x\bigr)^\dagger
A_\symm\bigl(\alpha^3 x+\tfrac{1-\alpha}{2}\bigr) L_3\,,
\equation(MatGFiG)
$$
and that the symmetric factor $\tilde A_\symm$ of $\tilde G$ is
$$
\eqalign{
\tilde A_\symm(x)
&=L_3^{-1}A_\symm\bigl(\alpha^3x+(1-\alpha)\bigr)^\dagger
B_\symm\bigl(\alpha^3x+\tfrac{1-\alpha}{2}\bigr)
A_\symm\bigl(\alpha^3x\bigr)^\dagger\times\cr
&\quad\times
B_\symm\bigl(\alpha^3x-\tfrac{1-\alpha}{2}\bigr)
A_\symm\bigl(\alpha^3x-(1-\alpha)\bigr)^\dagger L_3\,.\cr}
\equation(MatGiFGiFGi)
$$
Notice that, if $B_\symm$ and $A_\symm$
are analytic in a strip $|\Im x|<\delta$,
then $\tilde B_\symm$ and $\tilde A_\symm$
are analytic in the strip $|\Im x|<\alpha^{-3}\delta$.
The following is equally straightforward to verify.

Let $\radB$ and $\radA$ be positive real numbers.

\claim Proposition(ZerosOne)
Assume that $\alpha^3\radA+\thalf\alpha^2\le\radB
\le\alpha^{-3}\radA-\thalf\alpha^{-1}$.
Notice that $\radB=\sfrac{\alpha}{2}$
and $\radA=\shalf$ satisfy this condition.
If the domains of $B_\symm$ and $A_\symm$ include
the (real or complex) disks $|x|\le\radB$
and $|x|\le\radB$, respectively,
then the domains of $\tilde B_\symm$ and $\tilde A_\symm$
include the same disks.
Assume now that
$$
\radB>\tfrac{\alpha}{2}\,,\qquad
\radA>\thalf\,,\qquad
\alpha^3\radA+\thalf\alpha^2<\radB
<\alpha^{-3}\radA-\thalf\alpha^{-1}\,.
\equation(rhoDomainCond)
$$
Suppose that the domains of $B_\symm$ and $A_\symm$ include
the disks $|x|<\radB$ and $|x|<\radB$, respectively.
Then the domains of the symmetric factors
associated with $\RR_3^n(F,G)$ include the same disks.
Furthermore, these domains grow asymptotically
like $\alpha^{-3n}$, as $n\to\infty$.
Moreover, the values of these factors on any given compact set
are determined (for sufficiently large $n$) by the values
of $B_\symm$ on the disk $|x|<\radB$
and the values of $A_\symm$ on the disk $|x|<\radA$.

This motivates the following choice of function spaces.
Given $\rad>0$, denote by $\FF_\rad$ the space of all
real analytic functions $f$ on $(-\rad,\rad)$ that have a finite norm
$$
\|f\|_\rad=\sum_{n=0}^\infty\bigl|f^{(n)}(0)\bigr|\,{r^n\over n!}\,.
\equation(AANorm)
$$
Notice that every function $f\in\FF_\rad$
extends analytically to the complex disk $|x|<\rad$.
Furthermore, $\FF_\rad$ is a Banach algebra
under pointwise multiplication of functions.
This was crucial in [\rKochAM,\rKochUniv] but is less important here.

The space of matrix functions
$$
x\mapsto A_\symm(x)=\twomat{a_\symm(x)}{b_\symm(x)}{c_\symm(x)}{d_\symm(x)}
\equation(axbxcsxdx)
$$
with entries in $\FF_\rad$ will be denoted by $\FF_\rad^4$.
The norm of $A_\symm\in\FF_\rad^4$ is defined as
$\|A_\symm\|_\rad=\|a_\symm\|_\rad+\|b_\symm\|_\rad+\|c_\symm\|_\rad+\|d_\symm\|_\rad$.
Given a pair $\rad=(\radB ,\radA)$ of positive real numbers,
we define $\HH_\rad$ to be the vector space of all pairs
$\PP=(B_\symm,A_\symm)$ in $\FF_\radB^4\times\FF_\radA^4$,
equipped with the norm
$\|\PP\|_\rad=\|B_\symm\|_\radB+\|A_\symm\|_\radA$.

\smallskip
Given that every skew-product map that appears in our analysis
has a pre-determined first component,
we will identify a skew-product map $G=(\alpha,A)$
with its symmetric factor $A_\symm$.
Referring to \equ(axbxcsxdx),
we note that $G$ is reversible if and only if
$a_\symm$ and $d_\symm$ are even, while $b_\symm(x)=-c_\symm(-x)$.
The subspace of reversible pairs in $\HH_\rad$
will be denoted by $\HH_\rad'$.

\demo Convention(DomainCond)
We always assume that domain parameters $\rad=(\radB,\radA)$
satisfy the domain condition \equ(rhoDomainCond).
And whenever $L_k:\HH_\rad\to\pm\rmSL(2,\real)$ remains unspecified,
it is assumed to be continuous.

\claim Lemma(CompactRG)
The transformation $\RR_3:\HH_\rad\to\HH_\rad$ is compact
and preserves reversibility.

The reversibility-preserving property is a consequence of the fact that
the products $GF^\dagger G$ and $G^\dagger FG^\dagger FG^\dagger$
that appear in \equ(CCThree) are palindromic.
Compactness is a consequence of the fact that $\RR_3$
is analyticity-improving, in the sense that $\RR_3$ maps
bounded sets in $\HH_\rad$ to bounded sets in $\HH_{\rad'}$,
for some $\rad_\ssB'>\radB$ and $\rad_\ssA'>\radA$.
Now combine this with the fact that the inclusion maps
$\FF_{\rad_\ssB'}\to\FF_\radB$
and $\FF_{\rad_\ssA'}\to\FF_\radA$ are compact.

\smallskip
The same holds of course for the transformations $\RR_{3n}$ with $n>1$.

When the exact growth of products is not needed,
we also consider a modified version of our \RG transformation,
defined by
$$
\buR_{3n}=\buN\circ\RR_{3n}\,,
\equation(NormalizedRG)
$$
where $\buN$ is a suitable normalization.
Unless specified otherwise, $\buN$ consists
in normalizing $B_\symm\mapsto\|B_\symm\|_\radB^{-1}B_\symm$
and $A_\symm\mapsto\|A_\symm\|_\radA^{-1}A_\symm$.

The following implies \clm(subseqLimit).

\claim Theorem(GenScalingLimit)
Pick $L_3=S$. Let $P'=(F',G')$ be a pair in $\HH_\rad$,
of the form $F'=(1,\idmat)$ and $G'=(\alpha,A')$.
Assume hat $G'$ has a positive Lyapunov exponent and is reversible.
Denote by $K$ the set of all accumulation points
of the sequence $n\mapsto\buR_3^n(P')$.
Then every pair $P=(F,G)$ in $K$ has the following property.
$F=(1,B)$ and $G=(\alpha,A)$, with
$B_\symm(x)=b_\symm(x)\ttA^\dagger$ and $A_\symm(x)=a_\symm(x)\ttA$ for all $x$,
where $b_\symm$ and $a_\symm$ are nonzero even entire functions
and $\ttA$ is a constant $2\times 2$ matrix of rank $1$.

\proof
First, notice that $SC_\symm(x)S=C_\symm(-x)^\dagger$
for any reversible map $H=(\gamma,C)$.
So the conjugacy by $L_3=S$ in \equ(MatGFiG) and \equ(MatGiFGiFGi)
just replaces each factor $C(x)$ by $C(-x)^\dagger$.
In particular, the leftmost factor in \equ(MatGiFGiFGi)
becomes $A_\symm\bigl(-\alpha^3x-(1-\alpha)\bigr)$.

Since the transformation $\buR_3$ is compact, the set $K$ is compact and non-empty.
And by \clm(ZerosOne), the symmetric factors $B_\symm$ and $A_\symm$
of any pair $P\in K$ are entire analytic.
They have norm $1$ and thus cannot be zero.
Without the normalization $\buN$, the norm of $\RR_3^n(P_0)$
tends to infinity, since $G$ has a positive Lyapunov exponent.
So the factors $B_\symm$ and $A_\symm$ have determinant zero and rank $1$.

Consider the representation \equ(axbxcsxdx) for $A_\symm$.
By reversibility, the functions $a_\symm$ and $d_\symm$ are even, while $b_\symm(x)=-c_\symm(-x)$.
Consider first the case where one of $a_\symm,b_\symm,c_\symm,d_\symm$ is identically zero.
Then $a_\symm d_\symm=b_\symm c_\symm=0$, implying that $A_\symm=a_\symm\bigl[{1~0\atop 0~0}\bigr]$
or $A_\symm=d_\symm\bigl[{0~0\atop 0~1}\bigr]$.
Using now \equ(MatGFiG) and \equ(MatGiFGiFGi),
we see that $\tilde P=\buR_3(P)$ has symmetric factors
with the properties described in \clm(GenScalingLimit).

Assume now that none of $a_\symm,b_\symm,c_\symm,d_\symm$ is identically zero.
Consider the ratio $r_0=c_\symm/a_\symm$.
This is a meromorphic function describing
the angle between the range of $A_\symm$ and one of the coordinate axis.
Denote by $r_n$ the corresponding ratio
for the pair $P_n=\buR_3^n(P)$. By \equ(MatGiFGiFGi) we have
$$
r_1(x)=r_0\bigl(-\alpha^3 x-(1-\alpha)\bigr)\,.
\equation(tilderx)
$$
Notice that the function $x\mapsto-\alpha^3 x-(1-\alpha)$
has $x_0=-\sfrac{\alpha}{2}$ as a fixed point.
Setting $R_n(t)=r_n\bigl(x_0+t)$, the identity \equ(tilderx) becomes
$$
R_n(t)=R_0\bigl((-\alpha)^{3n} t\bigr)\,,
\equation(tildeRx)
$$
with $n=1$. The same holds for any $n>0$ by iteration.
Unless $R$ has a pole at the origin,
we see that $R_n\to R_0(0)$ on some open neighborhood of the origin.
If $R$ has a pole at the origin, we repeat the above with
$r_0=a_\symm/c_\symm$ in place of $r_0=c_\symm/a_\symm$.
Using compactness again,
this shows that $K$ includes a pair $P$ with constant
ratio $r_0$ in some open neighborhood of $x_0$.
In the case where $r_0=a_\symm/c_\symm$, this implies that
$a_\symm=r_0c_\symm$ and $b_\symm=r_0 d_\symm$ on all of $\complex$.
If $r_0=0$, then we are in the case discussed at the beginning.
Otherwise, the reversibility property $b_\symm=-c_\symm$
shows that $a_\symm,b_\symm,c_\symm,d_\symm$ are all constant multiples of $c_\symm$.
A similar argument applies in the case $r_0=c_\symm/a_\symm$.
In either case, since $a_\symm$ and $d_\symm$ are even
and not both identically zero,
$K$ includes a pair $P$ that admits the representation
described in \clm(GenScalingLimit).
\qed

We are interested in describing the points $x$
where the limiting factors $B_\symm$ and $A_\symm$
in \clm(GenScalingLimit) are zero.
For a limit pair $P=(F,G)$ with $F=(1,B)$ and $G=(\alpha,A)$,
the equations \equ(MatGFiG) and \equ(MatGiFGiFGi)
can be used to describe how zeros propagate
under renormalization.
Let $\BB$ be the set of zeros of $B_\symm$
and $\AA$ the set of zeros of $A _\symm$.
Let $\tilde\BB$ and $\tilde\AA$ be the set of zeros
of the symmetric factors for $\tilde P=\RR_3(P)$.
{}From \equ(MatGFiG) and \equ(MatGiFGiFGi) we see that
$$
\tilde\BB=\alpha^{-3}\hat\BB\,,\qquad
\hat\BB=\BB\bigcup_{m=\pm1}\Bigl(\AA+\tfrac{m}{2}\alpha^2\Bigr)\,,
\equation(LLtildeBBOne)
$$
and
$$
\tilde\AA=\alpha^{-3}\hat\AA\,,\qquad
\hat\AA=\bigcup_{m=\pm1}\Bigl(\BB+\tfrac{m}{2}\alpha^2\Bigr)
\bigcup_{m=0,\pm1}\Bigl(\AA+m\alpha^2\Bigr)\,.
\equation(LLtildeAAOne)
$$
The map $(\BB,\AA)\mapsto\bigl(\tilde\BB,\tilde\AA\bigr)$
will be described in more detail in Section 4.

\section Limit skew-products and the rotation number

In order to describe the behavior of factors \equ(StdAM)
with large values of $\lambda$,
we divide these matrices by $\lambda$ and consider the resulting factors
$$
A(x)=\twomat{v(t+x)-\epsilon}{-\delta}{\delta}{0}\,,
\equation(scaledSchrFac)
$$
where $\delta=\lambda^{-1}$ and $\epsilon=\lambda^{-1}E$.
Here, and in what follows, we assume that $\lambda>1$.
Orbits for the associated skew-product map $G=(\alpha,A)$ on $\torus\times\real^2$
are in one-to-one correspondence
with solutions of the equation $H_\delta^t u=\epsilon u$, where
$$
(H_\delta^t u)_q=v(t+q\alpha)u_q-\delta u_{q+1}-\delta u_{q-1}\,,\qquad q\in\integer\,.
\equation(scaledAMHam)
$$
Clearly, the operator $H_\delta^t$ on $\ell^2(\integer)$
has a well-defined limit as $\delta\to 0$,
given by setting $\delta=0$ in the above equation.
A sequence $u\in\ell^2(\integer)$ is an eigenvector of $H_0^t$
with energy $\epsilon$
if and only if $[v(t+q\alpha)-\epsilon]u_q=0$ for all $q$.
So for every integer $q$, the function $u=\delta_q$ supported at $q$
is an eigenvector of $H_0^t$ with eigenvalue $v(t+q\alpha)$.
Given that these eigenvectors span a dense subspace of $\ell^2(\integer)$,
the spectrum of $H_0^t$ is the closure of the range of $v$,
and it consists purely of eigenvalues and their accumulation points.

The integrated density of states $\buk(\epsilon)$
is the fraction (in the usual limit sense) of integers $q$
for which $v(t+q\alpha)\le\epsilon$.
By ergodicity, this is simply the measure of the set
of all points $x\in\torus$ for which $v(x)\le\epsilon$.

Let us compute $\buk(\epsilon)$ for the \AM potential
$v(x)=-2\cos(2\pi x)$. For simplicity consider $t=\sfrac{\alpha}{2}$.
Then the symmetric limit factor associated with \equ(scaledSchrFac) is
$$
A_\symm(x)=a_\symm(x)\stwomat{1}{0}{0}{0}\,,\qquad
a_\symm(x)=-\epsilon-2\cos(2\pi x)\,.
\equation(aoxLim)
$$
Given $\epsilon$ in the spectrum $[-2,2]$ of $H_0^t$,
we define $\rho=\rho(\epsilon)$ in $[0,\shalf]$ by the equation
$\epsilon=-2\cos(2\pi\rho)$.
Then
$$
a_\symm(x)=2\cos(2\pi\rho)-2\cos(2\pi x)=4\sin(\pi(x-\rho))\sin(\pi(x+\rho))\,.
\equation(acircmuepsilonTwo)
$$
If we represent $\torus$ as the interval $[-\shalf,\shalf]$
with the endpoints identified, then $a_\symm$ is negative on $(-\rho,\rho)$
and positive on the complement of $[-\rho,\rho]$.
This shows that $\buk(\epsilon)=2\rho$.

In order to determine the rotation number
of the skew-product map $G_\epsilon=(\alpha,A)$
with symmetric factor \equ(aoxLim),
let us first reformulate the definition of the number $\Rot_\ssN(G)$
that appears in \equ(rotviaSigns).
It counts half the number of times that the first
component of the vector $A^{\ast k}(t)y$ changes sign
as $k$ is increased from $1$ to $N$.
The limit \equ(rotviaSigns) is independent of $t$
and of the initial condition $y\in\real^2\setminus\{0\}$.
We choose the same definition for degenerate skew-products
with $A(x)=a_\symm(x)\ttA$, but we assume that $\ttA y\ne 0$.

Notice that $A^{\ast k}(t)y$ changes sign as $k$ is increased to $k+1$
precisely when $a_\symm(t+k\alpha)$ is negative.
So by ergodicity, the average number of sign changes
is the measure of $[-\rho,\rho]$.
The limit \equ(rotviaSigns) is half this number, which is $\rho$.
So in summary, we have
$$
\buk(\epsilon)=2\rho\,,\qquad
\rot(G_\epsilon)=\rho\,,\qquad \epsilon=-2\cos(2\pi\rho)\,.
\equation(trivIDSaandRot)
$$
In particular, $\buk(\epsilon)=2\rot(G_\epsilon)$,
which corresponds to the known relation for the \AM factors.
Notice that $H_0^t+2\delta\id\ge H_\delta^t\ge H_0^t-2\delta\id$,
where $H\ge K$ means that $H-K$ is a positive operator.
Using the definition of $\buk_\ssL$ given in the introduction, we have
$\buk_\ssL\bigl(H_0^t+2\delta\id,\epsilon\bigr)
\le\buk_\ssL\bigl(H_\delta^t,\epsilon\bigr)
\le\buk_\ssL\bigl(H_0^t-2\delta\id,\epsilon\bigr)$.
Taking $L\to\infty$ and then $\delta\to 0$,
and using that $\buk(\epsilon)=2\rot(G_\epsilon)$ for $0\le\delta<1$,
we obtain the following.

\claim Proposition(IDSConvergence)
Let $G_\delta=(\alpha,A)$, with $A$ given by \equ(scaledSchrFac)
and $v=-2\cos(2\pi\,\bdot)$.
Then $\rot(G_\delta)\to\rot(G_0)$ as $\delta\to 0$,
uniformly in $\epsilon$.

The same holds for real analytic potentials
close to $v=-2\cos(2\pi\,\bdot)$.

Notice that a skew-product map $G=(\alpha,A)$
with $A(x)=a(x)\bigl[{1~0\atop 0~0}\bigr]$ is for practical purposes
a skew-product map $g=(\alpha,a)$ on $\torus\times\real$,
where $a(x)$ acts on $\real$ by multiplication.
The rotation number of $g$, determined via sign changes,
is clearly the same as the rotation number of $G$.

\smallskip
Recall our definition \equ(rotviaSigns) of the rotation number $\rot(G)$
is restricted to Schr\"odinger factors as defined by \equ(StdAM).
The formula \equ(rotviaSigns) can be applied to other
pairs $G=(\alpha,A)$ as well.
But it does not yield the desired result, in general, not even modulo $\shalf$.

For a more general definition of a rotation number
for skew-product maps $G=(\alpha,A)$ on $\torus\times\real^2$,
assume that $A:\torus\to\rmSL(2,\real)$ is continuous.
Consider the map $y\mapsto A(x)y$ on $\real\setminus\{0\}$.
Writing $y$ as $\|y\|\bigl[{\cos t\atop\sin t}\bigr]$
and $\eta=A(x)y$ as $\|\eta\|\bigl[{\cos\tau\atop\sin\tau}\bigr]$,
we have a continuous function $\tau=\tau(x,t)$
from $\torus\times\circle$ to $\circle$,
where $\circle=\real/(2\pi\integer)$.
Each of the maps $t\mapsto\tau(x,t)$ on $\circle$ admits a continuous lift
$t\mapsto t+g(x,t)$ to $\real$.
This lift is unique up to an additive constant $2\pi m_x$ with $m_x\in\integer$.
Assuming now that $A$ is homotopic to $x\mapsto\idmat$,
$g$ can be chosen continuous, with $m_x=m$ for all $x$.
Using the map $\GG$ defined by $\GG(x,t)=(x+\alpha,t+g(x,t))$,
one now defines
$$
\varrho(G)=\lim_{N\to\infty}{\Sigma_\ssN(G)\over N}\,,\qquad
\Sigma_\ssN(G)\defeq{1\over 2\pi}\sum_{n=1}^N g\bigl(\GG^n(x,t)\bigr)\,.
\equation(varrhoG)
$$
Here, $\alpha$ can be any irrational number.
Using the unique ergodicity of $x\mapsto x+\alpha$ on $\torus$,
it is possible to prove [\rJoMo,\rHerman]
that the limit \equ(varrhoG) exists, is independent of $(x,t)$,
and that convergence is uniform in $(x,t)$.
Due to the freedom of choosing $m\in\integer$,
the rotation number $\varrho(G)$ is unique only modulo $1$.

\smallskip
$\Sigma_\ssN(G)$ measures the amount of rotation
on the interval $\{1,2,\ldots,N\}$, in units of $1$ per revolution.
For Schr\"odinger maps $G=(\alpha,A)$,
we assume that the lift $g$ has been chosen
in such a way that $g(0,0)\in[0,2\pi)$.
Then $\Sigma_\ssN(G)$ agrees within an error $\le 1$
with half the number of sign changes $\Rot_\ssN(G)$.
This follows from the fact that the second row of $A(x)$ is $[1\;0]$,
implying that $A(x)$ maps the right (left) half-plane to the upper (lower) half-plane.
So every full revolution is associated with
a single sign change positive $\to$ negative
and a single sign change negative $\to$ positive;
see also [\rDelSou,\rJohnson,\rPuigi].
This implies in particular that $\rot(G)=\varrho(G)$
for Schr\"odinger factors \equ(StdAM) with continuous potentials $v$.
(Replacing $A$ by $-A$ in this case
yields $\varrho\mapsto\varrho-\shalf$ modulo $1$,
while sign counting would yield $\rot\mapsto\shalf-\rot$.)

\section Periodic orbits

In this section we consider pairs of skew-product maps $p=(f,g)$
with $f=(1,b)$ and $g=(\alpha,a)$,
where $b$ and $a$ are even real analytic functions
whose values $a(x)$ and $b(x)$ act on $\real$ by multiplication.
The \RG transformations $\RR_{3n}$
extend naturally to such pairs, if we define $b^\dagger=b$ and $a^\dagger=a$.
Conjugacy has no effect for scalar factors,
so $\Lambda_k=\Lambda^k$ for all $k$.
The condition that $a$ and $b$ be even corresponds to reversibility.

Our goal here is to construct periodic orbits under renormalization
and to describe some of their properties.

\demo Convention(SymmFact)
Throughout this section,
the second component $A$ in a pair $(\alpha,A)$
representing a skew-product map denotes the symmetric factor.
So the composition rule here is
$(\beta,B)(\alpha,A)=(\alpha+\beta,C)$ with
$C(x)=B(x+\sfrac{\alpha}{2})A(x-\sfrac{\beta}{2})$.

First we need to establish some properties
of the map $(\BB,\AA)\mapsto\bigl(\tilde\BB,\tilde\AA\bigr)$
defined by the equations \equ(LLtildeBBOne) and \equ(LLtildeAAOne).
A crucial ingredient in the proof below is the well-known identity
$$
p_k-q_k\alpha=(-\alpha)^{k+1}\,,\qquad
q_k=p_{k+1}\,,
\equation(pkqkalpha)
$$
where $p_k$ denotes the $k$-th Fibonacci number.

\claim Lemma(ZerosTwo)
Let $\BB_0=\{\}$.
Given $\AA_0\subset\real$, define
$\BB_n=\tilde\BB_{n-1}$ and $\AA_n=\tilde\AA_{n-1}$ for $n=1,2,3,\ldots$
For every $r>0$ define $\AA_n(r)=\AA_n\cap[-r,r]$ and $\BB_n(r)=\BB_n\cap[-r,r]$.
Consider now $\AA_0=\{\rho\}$ with $0\le\rho\le\shalf$.
Then $\BB_n(\shalf)$ and $\AA_n(\shalf)$ contain at most one point, for each $n$.
The sequence $k\mapsto\AA_k(\shalf)$
is periodic if and only if $\rho$ is positive periodic
or belongs to $\{0,\sfrac{\alpha}{2},\shalf\}$.
Furthermore, some set $\AA_n(\shalf)$ includes a point $\tilde\rho\in\half\integer[\alpha]$
if and only if $\tilde\rho-\rho\in\integer[\alpha]$.

\proof
Let $B_0=\{\}$ and $A_0=\{\rho\}$,
considered as sets on the circle $\torus=\real/\integer$.
Let $B_n=\hat B_{n-1}$ and $A_n=\hat A_{n-1}$ for $n=1,2,3,\ldots$,
using the map $(B,A)\mapsto\bigl(\hat B,\hat A\bigr)$
defined by \equ(LLtildeBBOne) and \equ(LLtildeAAOne).
Consider the functions $b_0,a_0:\torus\to\real$
defined by $b_0(x)=1$ and $a_0(x)=\sin(\pi(x-\rho))$.
If we set $b_n=\hat b_{n-1}$ and $a_n=\hat a_{n-1}$ for $n=1,2,3,\ldots$, where
$$
\eqalign{
\hat b(x)
&=a\bigl(x-\tfrac{1-\alpha}{2}\bigr)b(x)
a\bigl(x+\tfrac{1-\alpha}{2}\bigr)\,,\cr
\hat a(x)
&=a\bigl(x+(1-\alpha)\bigr)
b\bigl(x+\tfrac{1-\alpha}{2}\bigr)a\bigl(x\bigr)
b\bigl(x-\tfrac{1-\alpha}{2}\bigr)
a\bigl(\alpha^3x-(1-\alpha)\bigr)\,,\cr}
\equation(idhataobo)
$$
then the zeros of $b_n$ and $a_n$
are precisely the sets $B_n$ and $A_n$, respectively.
Notice that
$$
\eqalign{
b_n(x)
&=a_0\bigl(x+\tfrac{q_{3n-1}}{2}\alpha\bigr)\cdots
a_0\bigl(x+\thalf\alpha\bigr)
a_0\bigl(x-\thalf\alpha\bigr)\cdots
a_0\bigl(x-\tfrac{q_{3n-1}}{2}\alpha\bigr)\,,\cr
a_n(x)
&=a_0\bigl(x+\tfrac{q_{3n}-1}{2}\alpha\bigr)\cdots
a_0\bigl(x+\alpha\bigr)a_0(x)
a_0\bigl(x-\alpha\bigr)\cdots
a_0\bigl(x-\tfrac{q_{3n}-1}{2}\alpha\bigr)\,.\cr}
\equation(Ancirc)
$$
So the zeros of $a_n$ constitute an orbit on $\torus$
under translation $x\mapsto x+\alpha\;(\mod1)$.
This orbit has length $q_{3n}$ and its center point is at $\rho$.
By the three-gap theorem [\rSos,\rSur,\rSwi],
the gaps between adjacent points in $A_n$
can only take the values $\alpha^{3n-2}$, $\alpha^{3n-1}$, and $\alpha^{3n}$.
In fact, the gap $\alpha^{3n-2}$ cannot occur,
since it gets closed in the last step of the orbit.
So the gaps in the scaled orbit $\AA_n=\alpha^{-3n}A_n$
are either $\alpha^{-1}$ or $1$.
Thus, $\AA_n(\shalf)$ contains at most one point.
Similarly, the gaps in $B_n$ are no shorter $\alpha^{3n-1}$,
so the scaled orbit $\BB_n=\alpha^{-3n}B_n$
has at most one point in $[-\shalf,\shalf]$.

\smallskip
Consider now the condition that $\AA_n$ contains
a given real number $\tilde\rho$.
This condition is equivalent to $\tilde a_n(\tilde\rho)=0$,
where $\tilde a_n(x)=a_n\bigl(\alpha^{3n}x\bigr)$.
{}From \equ(Ancirc) we see that $\tilde a_n(\tilde\rho)=0$ exactly when
$$
\bigl[(-\alpha)^{3n}-1\bigr]\tilde\rho-\nu_2+\nu_1\alpha
=\rho-\tilde\rho\,,\qquad
|\nu_1|\le{q_{3n}-1\over 2}\,,
\equation(AnCircZeroOne)
$$
for some $\nu\in\integer^2$.
Using that $(p_{3n-1}-1)-p_{3n}\alpha=(-\alpha)^{3n}-1$,
this can be written as
$$
\bigl[(p_{3n-1}-1)-p_{3n}\alpha\bigr]\tilde\rho-\nu_2+\nu_1\alpha
=\rho-\tilde\rho\,,
\qquad|\nu_1|\le{q_{3n}-1\over 2}\,.
\equation(AnCircZeroTwo)
$$

If $\tilde\rho=0$, then \equ(AnCircZeroTwo)
can be satisfied from some $n$ if and only if $\rho\in\integer[\alpha]$.
And a straightforward computation shows that
$\rho=0$ yields $\AA_n(\shalf)=\{0\}$ for all $n$.
Consider now an arbitrary $\tilde\rho\in\half\integer[\alpha]$.
In this case, the equation \equ(AnCircZeroTwo) can be written as
$$
\left[{p_{3n-1}-1\over 2}-{p_{3n}\over 2}\alpha\right]2\tilde\rho-\nu_2+\nu_1\alpha
=\rho-\tilde\rho\,,
\qquad|\nu_1|\le{q_{3n}-1\over 2}\,.
\equation(AnCircZeroThree)
$$
Notice that $p_{3n-1}-1$ and $p_{3n}$ are even.
So we must have $\rho\in\tilde\rho+\integer[\alpha]$.
For $\rho=\shalf$ and $\rho=\sfrac{\alpha}{2}$,
an explicit computation shows that $\AA_n(\shalf)=\{\pm\rho\}$ for all $n$.
We note also that \equ(AnCircZeroThree) has no solution
for $\tilde\rho=\shalf-\sfrac{\alpha}{2}$.

\smallskip
Consider now the periodicity condition $\tilde a_n(\rho)=0$,
which corresponds to setting $\tilde\rho=\rho$ in \equ(AnCircZeroTwo).
Clearly, we need $\rho\in\rational[\alpha]$
for \equ(AnCircZeroTwo) to have a solution in this case.
So write $\rho={w\over v}+{u\over v}\alpha$ for some $v>0$.
Without loss of generality,
we can restrict restrict to values of $n$
for which $p_{3n-1}\equiv 1$ and $p_{3n}\equiv 0$ modulo $v$.
This condition is satisfied e.g.~if $n$ is a multiple
of the Pisano period $\ell(v)$.
By definition, the Pisano period $\ell(v)$
is the period of the Fibonacci sequence $k\mapsto p_k$ modulo $v$.
Then we have $\tilde a_n(\rho)=0$, if and only if
$$
\left[{p_{3n-1}-1\over v}-{p_{3n}\over v}\alpha\right](w+u\alpha)-\nu_2+\nu_1\alpha=0\,,
\qquad|\nu_1|\le{q_{3n}-1\over 2}\,,
\equation(AnCircZeroThree)
$$
for some $\nu\in\integer^2$.
Multiplying out the product $[\ldots](w+u\alpha)$ in the above equation
and using that $\alpha^2=1-\alpha$, one finds that
$$
\nu_1={p_{3n}\over v}w-{q_{3n}-1\over v}u\,,\qquad
\nu_2={p_{3n-1}-1\over v}w-{p_{3n}\over v}u\,.
\equation(AnCircZeroThreenu)
$$
And the condition $|\nu_1|\le{q_{3n}-1\over 2}$ becomes
$$
\left|{p_{3n}\over q_{3n}-1}{w\over v}-{u\over v}\right|\le\half\,.
\equation(nuOneCond)
$$
Assume now that $\rho$ is neither $0$ nor $\shalf$,
since these values have been covered already.
Then equality in \equ(nuOneCond) can hold only for
finitely many values of $n$.
Given that ${p_{3n}\over q_{3n}-1}\to\alpha$,
we have $\tilde a_n(\rho)=0$ as $n\to\infty$
along multiples of $\ell(v)$
if and only if $\bigl|{w\over v}\alpha-{u\over v}\bigr|<\half$,
meaning that $\rho$ is positive periodic or equal to $\sfrac{\alpha}{2}$.
This show that the sequence $k\mapsto\AA_k(\shalf)$
is periodic if and only if $\rho$ is positive periodic
or belongs to $\{0,\sfrac{\alpha}{2},\shalf\}$.
\qed

We note that, since points evolve independently,
having $m>1$ points in $\AA_0$ results in $\AA_n(\shalf)$
and $\BB_n(\shalf)$ containing $m$ points or none.
Furthermore, if $\AA_0$ is invariant under $t\mapsto-t$,
then so are all the sets $\AA_n$ and $\BB_n$.

\demo Convention(fixrhon)
In the remaining part of this paper,
$\rho$ is a fixed but arbitrary positive periodic rotation number.
Referring to the periodicity statement in \clm(ZerosTwo),
the fundamental period of the sequence $k\mapsto\AA_k$,
generated from $\BB_0=\{\}$ and $\AA_0=\{\rho\}$,
will from now on be denoted by $n$.

Recall from \equ(NormalizedRG) that $\buR_{3n}=\buN\circ\RR_{3n}$,
where $\buN$ is a suitable multiplicative normalization.
The following theorem deals only with (zeros of) specific functions,
so we ``normalize'' here simply by identifying
functions that are constant multiples of each other.

\claim Theorem(PeriodicOne)
$\buR_{3n}$ has a fixed point $p_\diamond=((1,b_\diamond),(\alpha,a_\diamond)$
with the following properties.
The functions $b_\diamond$ and $a_\diamond$ are entire and even,
and as functions of $\sqrt{x}$, they are of order $\shalf$.
The zeros of $b_\diamond$ and $a_\diamond$ are all real and simple.
$a_\diamond$ has zeros at $\pm\rho$ and no other zeros
in the interval $[-\shalf,\shalf]$.

Following an approach used in [\rKochTrig],
we construct the functions $b_\diamond$ and $a_\diamond$
by first determining their zeros.
In what follows, we identify nonzero functions
$u$ and $v$ that are constant multiples of each other;
in symbols, $u\equiv v$.
The set of zeros of a function $u$ will be denoted by $\ZZ(u)$.

\smallskip
By reversibility, the factors of our skew-product maps
$f=(1,b)$ and $g=(\alpha,a)$ admit a factorization
$b(x)\equiv\bub(x)\bub(-x)$ and $a(x)\equiv\bua(x)\bua(-x)$.
This factorization is canonical (via orbits) and propagates under renormalization.
The functions $\bub$ and $\bua$ will be referred to as
the semi-factors of $f$ and $g$, respectively.
Consider the sequence of zeros $j\mapsto b_j$ of $\bub$
and the sequence of zeros $j\mapsto a_j$ of $\bua$.
If the sums $\sum_j|b_j|^{-2}$ and $\sum_j|a_j|^{-2}$ are finite, then
$$
b(x)\equiv\prod_j\biggl(1-{x^2\over b_j^2}\biggr)\,,\qquad
a(x)\equiv\prod_j\biggl(1-{x^2\over b_j^2}\biggr)\,,
\equation(ProdRep)
$$
with uniform convergence on compact subsets of $\complex$.
In the cases considered below, all zeros are simple,
so we may work with the sets $\BB=\ZZ(\bub)$ and $\AA=\ZZ(\bua)$
instead of the sequences $j\mapsto b_j$ and $j\mapsto a_j$, respectively.

If $\PP=(\BB,\AA)$ is a pair of subsets of $\real$,
define the pair $\tilde\PP=\bigl(\tilde\BB,\tilde\AA\bigr)$
by the equations \equ(LLtildeBBOne) and \equ(LLtildeAAOne).
Denoting the map $\PP\mapsto\tilde\PP$ by $R$,
we define $\PP_t=R^{tn}(\PP)$ for $t=0,1,2,\ldots$.
To recall the connection of $R$ with the transformation $\RR_3$,
consider two reversible maps $f=(1,b)$ and $g=(\alpha,a)$.
If $\BB=\ZZ(\bub)$ and $\AA=\ZZ(\bua)$, then we have
$\BB_t=\ZZ\bigl(\bub^t\bigr)$ and $\AA_t=\ZZ\bigl(\bua^t\bigr)$.
Here $\bub^t$ and $\bua^t$ are the semi-factors
associated with the pair $p^t=\RR_3^{tn}(p)$, where $p=(f,g)$.

\proofof(PeriodicOne)
In order to simplify the description,
consider the explicit pair $p=((1,1),(\alpha,a))$
with $a(x)=4\sin(\pi(x-\rho))\sin(\pi(x+\rho))$.
Its semi-factors are given by
$$
\bub(x)=1\,,\qquad\bua(x)\equiv\sin(\pi(x-\rho))\,.
\equation(AMab)
$$
Notice that $\bua$ is periodic with period $2$.
But its zeros are periodic with period $1$, that is,
$\BB=\{\}$ and $\AA=\{\rho\}+\integer$.
Since the transformation $R$ involves a scaling by $\alpha^{-3}$,
the sets $\BB_t$ and $\AA_t$ are periodic
with period $r_t=\alpha^{-3nt}$.
So we can regard them as subsets of the circle $\torus_t=\real/(r_t\integer)$.
Representing this circle by the interval $[-r_t/2,r_t/2]$
whose endpoints are being identified,
the zero sets $\BB_t'\subset\torus_t$ and $\AA_t'\subset\torus_t$
associated with $p^t$ are given by $\BB_t'=\BB_t\cap\torus_t$
and $\AA_t'=\AA_t\cap\torus_t$, respectively.

As described in the proof of \clm(ZerosTwo),
the set $\AA_t'$ is a scaled orbit of length $q_{3nk}$
under $x\mapsto x+\alpha$, centered at $x=\rho$.
As a result, the gaps between two adjacent zeros in $\AA_t'$
are either $\alpha^{-1}$ or $1$.
Similarly, the gaps between two adjacent zeros in $\BB_t'$
are no shorter than $\alpha^{-1}$.
Recall also that $\AA_0'\subset\AA_1'$ by \clm(ZerosTwo).
Thus, given that $\BB_0'\subset\BB_1'$, the inclusions
$\BB_t'\subset\BB_{t+1}'$ and $\AA_t'\subset\AA_{t+1}'$ hold for all $t$.

Now we can define $\BB_\infty=\lim_t\BB_t$ and $\AA_\infty=\lim_t\AA_t$,
using either the $\liminf$ or $\limsup$, with the same result.
Due to the above-mentioned gap properties,
the zeros $b^t_j$ of $\bub^t$ and the zeros $a^t_j$ of $\bua^t$
satisfy $\sum_j (b^t_j)^{-2}<\infty$ and $\sum_j(a^t_j)^{-2}<\infty$, respectively.
This includes the case $t=\infty$ as well.
We note that all zeros are real.
And they are simple: we have $\rho\not\in\half\integer[\alpha]$,
so $-\rho$ is not on the orbit of $\rho$
under translation $x\mapsto x+\alpha$ on $\torus$.

Clearly, $\PP_\infty=(\BB_\infty,\AA_\infty)$ is a fixed point of $R^{3n}$.
Now define $b_\diamond$ and $a_\diamond$ as products \equ(ProdRep),
for the zeros $b_j=b^\infty_j$ and $a_j=a^\infty_j$, respectively.
Since we are dealing with a class of entire functions
\big(even, and of order $\shalf$ in the variable $\sqrt{x}$ \big)
that are determined up to a constant factor by their zeros,
it is not hard to see that the pair $p_\diamond$ with symmetric factors
$b_\diamond$ and $a_\diamond$ is a fixed point of $\RR_{3n}$,
up to a normalization by constant factors.
\qed

We note that $b^t\to b_\diamond$ and $a^t\to a_\diamond$,
uniformly on compact subsets of $\complex$.
This is a straightforward consequence of our uniform gap-bounds,
together with the property that
$\BB_t'\subset\BB_{t+1}'$ and $\AA_t'\subset\AA_{t+1}'$ for all $t$.
But we need not prove this here,
since it follows from Theorem 4.3 below.

For a more detailed description of the convergence argument
we refer to [\rKochTrig], where a similar analysis has been carried out
(related to the critical case)
for meromorphic factors $a(x)=\bma(x)/\bma(-x)$ and quadratic irrational $\alpha$.

We note that the product \equ(ProdRep)
for the function $a_\diamond$ is in effect an infinite product
of sine functions, albeit rescaled.
The classical product of sines is the Sudler product
$P_m(\alpha)=\prod_{j=1}^m|2\sin\pi j\alpha|$.
Perturbed versions of this product, for $\alpha$ the golden mean,
have been investigated recently in [\rATZ].

\smallskip
Let us return to the issue of normalization.
As can be seen from \equ(MatGFiG) and \equ(MatGiFGiFGi),
we may replace $b_\diamond\mapsto-b_\diamond$
and/or $a_\diamond\mapsto-a_\diamond$,
in such a way that $b_\diamond(0)>0$ and $a_\diamond(0)<0$.
This corresponds to the signs of the limit \AM pair,
with $b_\symm=1$ and $a_\symm$ given by \equ(aoxLim).

The symmetric factors of $\tilde p_\diamond=\RR_{3n}(p_\diamond)$
are of the form $\tilde b_\diamond=b_\diamond e^{\tilde v_\diamond}$
and $\tilde a_\diamond=a_\diamond e^{\tilde u_\diamond}$
for some real constants $\tilde v_\diamond$ and $\tilde u_\diamond$.
Consider a pair $p=((1,a),\alpha,b))$
with symmetric factors $b=b_\diamond e^v$ and $a=a_\diamond e^u$.
If $u$ and $v$ are constants,
then the symmetric factors for $\tilde p=\RR_{3n}(p)$
are $\tilde b=\tilde b_\diamond e^{\tilde v}$
and $\tilde a=\tilde a_\diamond e^{\tilde u}$, where
$$
\twovec{\tilde v}{\tilde u}
=\twovec{\tilde v_\diamond}{\tilde u_\diamond}
+U^{3n}\twovec{v}{u}\,,\qquad U=\twomat{0}{1}{1}{1}\,.
\equation(expRG)
$$
Since $U$ is hyperbolic, it is clear that
we can determine $v$ and $u$ in such a way that $\tilde v=v$ and $\tilde u=u$.
Thus, we will assume from now on that $b_\diamond$ and $a_\diamond$
have been normalized in such a way that $p_\diamond$
is a fixed point of $\RR_{3n}$.

For the normalization $\buN$ that enters
the definition \equ(NormalizedRG) of $\buR_{3n}$,
we now choose $(b,a)\mapsto\bigl(Mb,a\bigr)$,
with $M$ defined by the condition $Mb(0)=b_\diamond(0)$.
Clearly, $p_\diamond$ is a fixed point of $\buR_{3n}$.

\smallskip
Denote by $\kappa$ the positive zero of $b_\diamond$
that is closest to the origin.
Presumably $\kappa>\sfrac{\alpha}{2}$, but we are not sure.
In any case, all zeros of $p_\diamond$ are determined (via orbits)
from the two zeros of $a_\diamond$ at $\pm\rho$.
Under renormalization, zeros of $b_\diamond$ in $[-\sfrac{\alpha}{2},\sfrac{\alpha}{2}]$
could only generate zeros of $b_\diamond$ outside
$[-\sfrac{\alpha}{2},\sfrac{\alpha}{2}]$
and/or zeros of $a$ outside $[\shalf,\shalf]$.

The following holds for some open neighborhood $Z$ of $\kappa$.

\claim Theorem(PeriodicOnePlus)
Let $b$ be even and real analytic on $[-\sfrac{\alpha}{2},\sfrac{\alpha}{2}]$,
satisfying $b(0)>0$.
Let $a$ be even and real analytic on $[-\shalf,\shalf]$,
satisfying $a(0)<0$.
Assume that $b$ has no zero in
$(0,\sfrac{\alpha}{2}]\setminus Z$.
Assume that $a$ has a simple zero at $\rho$
and no other zeros in $(0,\shalf]$.
If $p=((1,b),(\alpha,a))$,
then $\buR_{3n}^k(p)\to p_\diamond$ as $k\to\infty$,
uniformly on compact subsets of $\complex$.

Given any real numbers $\radA$ and $\radB$ satisfying \equ(rhoDomainCond),
define $\GG_\rad=\FF_\radB\times\FF_\radA$
and denote by $\GG_\rad'$ the subspace of reversible pairs.

\proofof(PeriodicOnePlus)
Consider first the sequence of pairs $p_m=\RR_{3n}^m(p)$.
For sufficiently large $m$,
the zero sets of $b_m$ and $a_m$
agree with the zero sets of $b_\diamond$ and $a_\diamond$, respectively,
on any given compact subset of $\complex$.
This includes the absence of non-real zeros, since their imaginary parts
expand by a factor $\alpha^{-3}$ under $\RR_3$,
as can be seen from \equ(LLtildeBBOne) and \equ(LLtildeAAOne).

Now pick $\radA$ and $\radB$ satisfying \equ(rhoDomainCond).
If $m_0$ has been chosen sufficiently large, then for $m\ge m_0$,
the pair $p_m$ belongs to $\GG_\rad'$,
and its zeros agree with those of $p_\diamond$
on the disks $|x|\le\radB$ for $b$ and $|x|\le\radA$ for $a$.
To simplify notation, replace now $p$ by $p_{m_0}$.

Then we have $b=b_\diamond e^{v}$ and $a=a_\diamond e^u$
for a pair of functions $(u,v)\in\GG_\rad$.
And the symmetric factors
of $\tilde p=\buR_{3n}(p)$ are again of the form
$\tilde b=b_\diamond e^{\tilde v}$ and $\tilde a=a_\diamond e^{\tilde u}$,
with $(\tilde u,\tilde v)\in\GG_\rad$.
The map $L: (v,u)\mapsto(\tilde v,\tilde u)$ on $\GG_r'$
is linear, compact, and rather simple:
its eigenvectors are pairs of polynomials of degree $2k$
with eigenvalues $\alpha^{3n(2k\pm 1)}$;
except at degree $0$, where our normalization condition $\buN$ acts.
Without $\buN$, the constant functions
transform via the matrix $U^{3n}$ in \equ(expRG).
Our normalization condition $\buN$
eliminates the expanding eigenvalue $\alpha^{-3n}$
and makes the constant terms contract as well.
Thus, $u$ and $v$ tend to zero under iteration of $L$,
implying that $\buR_{3n}^k(p)\to p_\diamond$ in $\GG_r'$,
as $k\to\infty$.
\qed

\demo Remark(NonPeriodic)
\clm(subseqLimit) and \clm(ZerosTwo) indicate that
the periodic orbits $p_\diamond$
associated with periodic rotation numbers $\rho$
are merely periodic points in a supercritical attractor for $\buR_3$.
It should be possible to determine non-periodic orbits
associated with many other values of $\rho$.
In such an analysis, there would be no reason to restrict $\alpha$
to the inverse golden mean or another quadratic irrational.
But we have not looked into these questions.

\section Hyperbolicity and consequences

Let $p_\diamond$ be the fixed point of $\buR_{3n}$
that was described in the last section.

\clm(PeriodicOnePlus) implies that the stable subspace
of $D\buR_{3n}(p_\diamond)$ has codimension one
and includes all reversible pairs
$p=((1,b),(\alpha,a))$ that satisfy $a(\rho)=0$.
Perturbations with $a(\rho)\ne 0$ are unstable:
\equ(LLtildeBBOne) and \equ(LLtildeAAOne)
show that small perturbations of zeros grow at a rate $\alpha^{-3}$
under $\RR_3$.
This implies the following.

\claim Corollary(PeriodicOneCor)
The derivative of $D\buR_{3n}(p_\diamond)$ on $\GG_r'$
has a simple eigenvalue $\alpha^{-3n}$
and no other spectrum outside the open unit disk.
The stable subspace is characterized by the condition $a_\symm(\rho)=0$.

Consider now the transformation $\RR_3$ for skew-products
with factors in $\rmGL(2,\real)$.
Given the form of our matrices \equ(scaledSchrFac),
we associate with the fixed point $p_\diamond$
the pair $P_\diamond=(F_\diamond,G_\diamond)$
with symmetric factors
$$
B_\diamond(x)=b_\diamond(x)\ttA^\dagger\,,\qquad
A_\diamond(x)=a_\diamond(x)\ttA\,,\qquad
\ttA=\stwomat{1}{0}{0}{0}\,.
\equation(Pdiamond)
$$
Recall from \equ(PaliRG) that the transformation $\RR_{3n}$
includes a scaling $\Lambda_{3n}(x,y)=\bigl(\alpha^{3n}x,L_{3n}y\bigr)$.
Consider first the choice $L_{3n}=S^n$.
Then it is straightforward to check that $P_\diamond$ is a fixed point of $\RR_{3n}$.
If we restrict $\RR_{3n}$ to pairs
$P=(F,G)$ with symmetric factors
$B_\symm(x)=b_\symm(x)\ttA^\dagger$ and $A_\symm(x)=a_\symm(x)\ttA$,
where $\ttA$ is fixed as in \equ(Pdiamond),
then $D\RR_{3n}(P_\diamond)$ is clearly equivalent to $D\RR_{3n}(p_\diamond)$.
But in the full space $\HH_\rad'$, the derivative of $D\RR_{3n}(P_\diamond)$
has an eigenvalue $1$
associated with conjugacies by constant matrices that commute with $S$.
This eigenvalue can be eliminated as follows.

Let $\bigl(\check F,\check G\bigr)=\RR_{3n}(F,G)$, still with the choice $L_{3n}=S^n$.
As an extra normalization step,
we include a conjugacy
$\check A_\symm\mapsto\tilde A_\symm=e^{-\sigma_{3n} S}\check A_\symm e^{\sigma_{3n} S}$
and
$\check B_\symm\mapsto\tilde B_\symm=e^{-\sigma_{3n} S}\check B_\symm e^{\sigma_{3n} S}$,
in order to normalize the off-diagonal terms $\tilde b_\symm$
and $\tilde c_\symm$ of the matrix $\tilde A_\symm$.
To be more precise, we choose $\sigma_{3n}=\sigma_{3n}(P)$
in such way that
$$
\tilde c_\symm(x_0)=0\,,\qquad x_0=-\tfrac{\alpha}{2}\,.
\equation(sigmaCond)
$$
This is possible if $\check A_\symm$ is sufficiently close to $A_\diamond$,
since $a_\diamond(x_0)\ne 0$.
Here we use that $\rho\ne x_0$.
Our reason for normalizing at $x_0$ is that this is the point
about which we have the scaling \equ(tildeRx).
This normalization step can now be incorporated into the scaling $\Lambda_{3n}$
by choosing
$$
L_{3n}=S^ne^{\sigma_{3n} S}=S^n\bigl[\cosh(\sigma_{3n})\idmat+\sinh(\sigma_{3n})S\bigr]\,,\qquad
\sigma_{3n}=\sigma_{3n}(P)\,.
\equation(newLii)
$$
In what follows, the transformation $\RR_{3n}$
is defined with the above choice of $L_{3n}$.

We also have to define a proper version of $\buR_{3n}=\buN\circ\RR_{3n}$.
This can be done the same way as for scalar factors.
Let $\tilde P=\buR_{3n}(P)$.
If $\tilde B_\symm$ denotes the symmetric factor of $\tilde F$,
then $\buN$ consists of the scaling $\tilde B_\symm\mapsto M\tilde B_\symm$,
where $M$ is determined by the condition that
$$
M\tr\bigl(\tilde B_\symm(0)\bigr)=b_\diamond(0)\,.
\equation(newNN)
$$
Clearly, $P_\diamond$ is a fixed point of $\buR_{3n}$.

Consider now the derivative of $\buR_{3n}(P_\diamond)$,
applied to a pair $P=(F,G)$ with symmetric factors
$$
B_\symm=\stwomat{a_\ssB}{b_\ssB}{c_\ssB}{d_\ssB}\,,\qquad
A_\symm=\stwomat{a_\ssA}{b_\ssA}{c_\ssA}{d_\ssA}\,.
\equation(tangBoAo)
$$

\claim Theorem(Hyperbolicity)
$D\buR_{3n}(P_\diamond)$ is hyperbolic,
with a simple eigenvalue $\alpha^{-3n}$
and no other spectrum outside the open unit disk.
The stable subspace of $D\buR_{3n}(P_\diamond)$
is characterized by the condition $a_\ssA(\rho)=0$
in the representation \equ(tangBoAo).

\proof  
Consider a representation analogous to \equ(tangBoAo)
for the pair $\tilde P=D\buR_{3n}(P_\diamond)P$,
with a tilde over the matrices and their elements.

The factors $\tilde B_\symm$ and $\tilde A_\symm$
naturally decompose into sums of products of matrices
from $\{B_\symm,A_\symm,B_\diamond,A_\diamond\}$.
But only a few products are nonzero, due to the fact
that each product can only have one factor $B_\symm$ or $A_\symm$.
The other factors all are $B_\diamond$ or $A_\diamond$,
which are ``sparse'' in the sense that they only have one nonzero matrix element.

In order to describe some nonzero terms,
let us split $B_\symm$ into a ``\onedim'' part that has $a_\ssB=b_\ssB=c_\ssB=0$,
and into a ``\twodim'' part that has $d_\ssB=0$.
Similarly, we split $A_\symm$ into a \onedim~ part
that has $b_\ssA=c_\ssA=d_\ssA=0$,
and into a \twodim~ part that has $a_\ssA=0$.
If we split $\tilde B_\symm$ and $\tilde A_\symm$ analogously,
then $D\buR_{3n}(P_\diamond)$ becomes a $2\times 2$ ``matrix''.
A useful feature of this matrix is that it is upper triangular,
meaning that the entry ``$\onedim\mapsto\twodim$'' is zero.
So the eigenvalues of $D\buR_{3n}(P_\diamond)P$
are those of the operators ``$\onedim\mapsto\onedim$''
and ``$\twodim\mapsto\twodim$''.
The former is hyperbolic, with a single expanding direction,
as described in \clm(PeriodicOneCor).

Consider now the ``$\twodim\mapsto\twodim$'' part.
The first observation is that $\tilde a_\ssB$ and $\tilde d_\ssA$
are zero when $\sigma_{3n}=0$.
So these functions are determined by $c_\ssA(x_0)$ only.
As a result, $D\buR_{3n}(P_\diamond)$ has an an eigenvalue $0$
with eigenvectors whose only nonzero components are $a_\ssB$ and $d_\ssA$.
Furthermore, the terms $a_\ssB$ and $d_\ssA$ do not contribute
to the off-diagonal entries of $\tilde B_\symm$ or $\tilde A_\symm$.
So it suffices to consider the operator ``off-diagonal $\mapsto$ off-diagonal''.
Here, only the rightmost or leftmost factors contribute,
and this is always a factor $A_\symm$ or $A_\symm^\dagger$.
So $\tilde c_\ssB$ and $\tilde b_\ssB$ are determined by $c_\ssA$ and $b_\ssA$.
And $D\buR_{3n}(P_\diamond)$ has an an eigenvalue $0$
with eigenvectors whose only nonzero components are $c_\ssB$ and $b_\ssB$.

Consider now the case $\sigma_{3n}=0$, where the conjugacy by $e^{\sigma_{3n} S}$ is trivial.
Then $\tilde c_\ssA$ satisfies
$$
{\tilde c_\ssA\over a_\diamond}(x_0+t)
={c_\ssA\over a_\diamond}\bigl(x_0-\alpha^{3n}t\bigr)\,,
\equation(dertildeRx).
$$
near $t=0$, where $x_0=-\sfrac{\alpha}{2}$.
This corresponds of course to the relation \equ(tildeRx).
So $D\buR_{3n}(P_\diamond)$ has eigenvalues $(-\alpha)^{3nm}$
with eigenvectors whose only nonzero components are
$$
c_\ssA(x)=a_\diamond(x)(x-x_0)^m\,,\qquad
c_\ssB(x)=-b_\diamond(x)(x-x_0)^m\,,
\equation(cAcBEigenvec)
$$
as well as $b_\ssA(x)=-c_\ssA(-x)$ and $b_\ssB(x)=-c_\ssB(-x)$.
These eigenfunctions are contracted for $m>0$.
But for $m=0$ we have an eigenvalue $1$.
This is due to the fact that we considered $\sigma_{3n}=0$.
With $\sigma_{3n}=\sigma_{3n}(P)$ as defined by the condition \equ(sigmaCond),
the eigenvalue becomes $0$.

This shows that the only unstable direction
of $D\buR_{3n}(P_\diamond)$ is the one inherited from $D\buR_{3n}(p_\diamond)$,
which is associated with a nonzero value of $a_\ssA(\rho)$.
So the stable subspace $W^s$ of $D\buR_{3n}(P_\diamond)$ consists
of all pairs $P$ with the property that $a_\ssA(\rho)=0$.
\qed

Our proof of \clm(supercritAMLimits) below involves
the pre-limit versions
$$
\rot_\ssN(G)={\Rot_N(G)\over N}\,,\qquad
\varrho_\ssN(G)={\Sigma_N(G)\over N}\,,
\equation(rotNvarrhoN)
$$
of the rotation numbers
defined in \equ(rotviaSigns) and \equ(varrhoG), respectively.
Let $P=(F,G)$ with $F=(1,\idmat)$.
Then the rotation numbers \equ(rotNvarrhoN) for $N=q_{3nk}$
can be computed via iterates $P_k=\buR_{3n}^k(P)$,
using that for even $k$,
$$
P_k=(F_k,G_k)\,,\qquad
G_k=\Lambda_{3nk}^{-1}\bigl(F^\dagger\bigr)^{p_{3nk}}G^{q_{3nk}}\Lambda_{3nk}\,.
\equation(GkThree)
$$
Here, $\Lambda_{3nk}(x,y)=\bigl(\alpha^{3nk}x,e^{\sigma_{3nk}S}y\bigr)$,
and $\sigma_{3nk}$ is the sum $\sigma_{3nk}(P)=\sum_{j<k}\sigma_{3n}(P_j)$.

Assuming that $G=(\alpha,A)$ commutes with $F$,
meaning that $A$ is periodic with period $1$,
the factors $F^\dagger$ in \equ(GkThree) contribute nothing
to the matrix part of $G_k$.
We will be using \equ(GkThree) in the case where $P_k\to P_\diamond$,
with the convergence being asymptotically geometric.
Then the sequence of normalization constant $k\mapsto\sigma_{3nk}(P)$ converges,
since $\sigma_{3n}(P_\diamond)=0$.

\proofof(supercritAMLimits)
Consider first the claim that \equ(rhoSymmFacSeq) holds
with $k\mapsto e^{-p_{3nk}L}$ and $k\mapsto e^{-q_{3nk}L}$
replaced by suitable sequences $k\mapsto M_k$ and $k\mapsto M_k'$, respectively.

Let $R_d$ be a small rectangle in $\HH_\rad'$, centered at $P_\diamond$.
To be more precise, let us call the stable (codimension $1$)
subspace $W^s$ of $D\buR_{3n}(P_\diamond)$ vertical
and the unstable (dimension $1$) subspace $W^u$ horizontal.
Then $R_d$ is a rectangle centered at $P_\diamond$,
of width and height $2d$, whose left/right sides are vertical
and top/bottom sides are horizontal.
Denote by $\WW^s$ ($\WW^u$) the local (un)stable manifold of $\buR_{3n}$ at $P_\diamond$.
It is tangent to $P_\diamond+W^s$ ($P_\diamond+W^u$) at $P_\diamond$.
If $d>0$ is chosen sufficiently small,
then $\WW^s$ ($\WW^u$) leaves the rectangle through the top/bottom (left/right)
sides within $\OO\bigl(d^2\bigr)$ from their centers.

Consider first the \AM factor,
multiplied by $\delta=\lambda^{-1}$, with energy $E=\lambda\epsilon$.
So we have a pair $P=(F,G)$ with $F=(1,\idmat)$ and $G=(\alpha,A)$,
and the symmetric factor $A_\symm$ is
$$
A_\symm(x)=\twomat{a_\symm(x)}{-\delta}{\delta}{0}\,,\qquad
a_\symm(x)=-\epsilon-2\cos(2\pi x)\,.
\equation(ScaledAM)
$$
To indicate the dependence on $\delta$ and $\epsilon$,
we will use the notation $P^\delta$ or $P^\delta(\epsilon)$.

Consider first the case $\delta=0$.
Let $\epsilon(0)=-2\cos(2\pi\rho)$.
If we fix $\epsilon=\epsilon(0)$,
then the sequence $k\mapsto\buR_{3n}^k\bigl(P^0\bigr)$ converges to
$P_\diamond$ by \clm(PeriodicOnePlus).
If we consider the family $\epsilon\mapsto P^0(\epsilon)$,
with $\epsilon$ close to $\epsilon(0)$,
then for sufficiently large $k$,
the family $\epsilon\mapsto\buR_{3n}^k(P^0(\epsilon))$
enters and leaves $R_d$ through the vertical sides
and it is transversal to the stable manifold $\WW^s$.
This follows from the way the zeros of $a_\symm$
transform under renormalization;
and we assume that $d>0$ has been chosen sufficiently small.

The same holds for the family $\epsilon\mapsto P^\delta(\epsilon)$,
for $\delta\ne 0$ sufficiently close to $0$.
So if $d>0$ has been chosen sufficiently small,
then the renormalized family $\epsilon\mapsto\buR_{3n}^k(P^\delta(\epsilon))$
crosses $\WW^s$ transversally at some value $\epsilon=\epsilon(\delta)$.
Clearly, the curve $\delta\mapsto\epsilon(\delta)$ is real analytic
near the origin, and it does not depend on $k$.

What remains to be shown is that
$G^\delta(\epsilon(\delta))$ has rotation number $\rho$.
To this end, let $P^\delta=P^\delta(\epsilon(\delta))$
and $P_k^\delta=\buR_{3n}^k\bigl(P^\delta\bigr)$.
To simplify notation, let us restrict to $\delta>0$.

\smallskip
We note that, if $k$ is sufficiently large,
then the symmetric factor $A^\delta_{k,\symm}(0)$ of $G^\delta_k$ is close to
$A_\diamond(0)=a_\diamond(0)\bigl[{1~0\atop 0~0}\bigr]$.
Thus, given that $a_\diamond(0)\ne 0$,
the first component of $A^\delta_{k,\symm}(0)\bigl[{1\atop 0}\bigr]$
stays bounded away from $0$ and thus does not change sign, as $\delta$ is varied.
This holds uniformly in $k$ and $\delta$,
if $k\ge k_0$ and $\delta\le\delta_0$ for some $k_0>0$ and some $\delta_0>0$.

We need a slight variation of this property.
Recall from \equ(GkThree) that $A^\delta_{k,\symm}(0)$ is
the product $\bigl(A^\delta\bigr)^{\ast q_{3nk}}_\symm(0)$,
conjugated by a matrix $e^{\sigma_{3nk}S}$
with $\sigma_{3nk}=\sigma_{3nk}\bigl(P^\delta\bigr)$.
Here, and in what follows, we assume that $k$ is even.
Notice that $\sigma_{3nk}\bigl(P^\delta\bigr)\to 0$ as $\delta\to 0$,
uniformly in $k$.
Thus, by decreasing the value $\delta_0>0$, if necessary,
the absence of sign changes described above
still holds if $A^\delta_{k,\symm}(0)$
is replaced by $\bigl(A^\delta\bigr)^{\ast q_{3nk}}_\symm(0)$.

Now choose a fixed but arbitrary positive $\delta\le\delta_0$.
Then the above implies that
$$
\rot_{N_k}\bigl(G^\delta\bigr)=\rot_{N_k}\bigl(G^0\bigr)\,,\qquad
N_k=q_{3nk}\,,\quad k\ge k_0\,.
\equation(ConstRot)
$$
Here, and in what follows,
we use as starting point $x={1-N\over 2}\alpha$
in the definition of $\Rot_\ssN(G)$.
This corresponds to the argument of the rightmost factor
in our symmetric products $A^{\ast\ssN}_\symm(0)$.
Due to the uniform convergence in \equ(varrhoG),
we still have $\Rot_\ssN(G)\to\rot(G)$ as $N\to\infty$.

Let now $\eps>0$.
If $k\ge k_0$ is sufficiently large, then
$\bigl|\rot_{N_k}\bigl(G^\delta\bigr)-\rot\bigl(G^\delta\bigr)\bigr|<\eps/2$
and $\bigl|\rot_{N_k}\bigl(G^0\bigr)-\rho\bigr|<\eps/2$.
Combining this with \equ(ConstRot), we have
$$
\eqalign{
\bigr|\rot\bigl(G^\delta\bigr)-\rho\bigl|
&\le\bigl|\rot\bigl(G^\delta\bigr)-\rot_{N_k}\bigl(G^\delta\bigr)\bigr|\cr
&\quad+\bigl|\rot_{N_k}\bigl(G^\delta\bigr)-\rot_{N_k}\bigl(G^0\bigr)\bigr|
+\bigl|\rot_{N_k}\bigl(G^0\bigr)-\rho\bigr|<\eps\,.\cr}
\equation(rhoBoundGdelta)
$$
Since $\eps>0$ was arbitrary, we conclude that $\rot\bigl(G^\delta\bigr)=\rho$.

\smallskip
Given a positive value $\delta_1\le\delta_0$,
the above can be extended to more general real analytic families of maps.
We have to assume that,
after a suitable rescaling and possible reparametrization,
the resulting family of factors $A^{\delta,\epsilon}_\circ$
is sufficiently close in $\FF_r^4$ to the rescaled \AM family \equ(ScaledAM),
uniformly in $\delta$ (close to $\delta_1$)
and $\epsilon$ (close to $\epsilon(\delta_1)$).
Then we can repeat the above argument for this family
(using $\Sigma_\ssN$ in place of $\Rot_\ssN$)
and get the same conclusion.
The main point here is that transversality to $\WW^s$
is stable under small perturbations.
\qed


\bigskip
\references

{\ninepoint

\item{[\rHarp]} P.G.~Harper,
{\sl Single band motion of conduction electrons in a uniform magnetic field},
Proc. Phys. Soc. Lond. A {\bf 68}, 874--892 (1955).

\item{[\rSos]} V.T.~S\'os,
{\sl On the distribution {\rm mod $1$} of the sequence $n\alpha$},
Ann. Univ. Sci. Budapest,
E\"otv\"os Sect. Math. {\bf 1}, 127--134 (1958).

\item{[\rSur]} J.~Sur\'anyi,
{\sl \"Uber die Anordnung der Vielfachen einer reelen Zahl {\rm mod $1$}},
Ann. Univ. Sci. Budapest,
E\"otv\"os Sect. Math. {\bf 1}, 107--111 (1958).

\item{[\rSwi]} S.~\'Swierczkowski,
{\sl On successive settings of an arc on the circumference of a circle},
Fundamenta Mathematicae {\bf 46}, 187--189 (1959).

\item{[\rKato]} T.~Kato,
{\sl Perturbation theory for linear pperators},
Springer Verlag, 1976.

\item{[\rHof]} D.R.~Hofstadter,
{\sl Energy levels and wave functions of Bloch electrons
in rational and irrational magnetic fields},
Phys. Rev. B {\bf 14}, 2239--2249 (1976).

\item{[\rTKNdN]} D.J.~Thouless, M.~Kohmoto, M.P.~Nightingale, M.~den Nijs,
{\sl Quantized Hall conductance in a two-dimensional periodic potential},
Phys. Rev. Lett. {\bf 49}, 405--408 (1982).

\item{[\rJoMo]} R.~Johnson, J.~Moser,
{\sl The rotation number for almost periodic potentials},
Commun. Math. Phys. {\bf 84}, 403--438 (1982).

\item{[\rHerman]} M.R.~Herman,
{\sl Une m\'ethode pour minorer les exposants de Lyapounov et quelques exemples
montrant le caract\`ere local d'un th\'eor\`eme d'Arnold et de Moser
sur le tore de dimension 2},
Comment. Math. Helvetici {\bf 58}, 453--502 (1983).

\item{[\rDelSou]} F.~Delyon, B.~Souillard,
{\sl The rotation number for finite difference operators and its properties},
Commun. Math. Phys. {\bf 89}, 415--426 (1983).

\item{[\rAvSi]} J.~Avron, B.~Simon,
{\sl Almost periodic Schr\"odinger operators. II.
The integrated density of states},
Duke Math. J. {\bf 50}, 369--391 (1983).

\item{[\rJohnson]} R.A.~Johnson,
{\sl A review of recent work on almost periodic differential and difference operator},
Acta Appl. Math. {\bf 1}, 241--261 (1983)

\item{[\rWieZa]} P.B.~Wiegmann, A.V.~Zabrodin,
{\sl Quantum group and magnetic translations.
Bethe ansatz solution for the Harper's equation},
Modern Phys. Lett. B {\bf 8}, 311--318 (1994).

\item{[\rBeGap]} J.~Bellissard,
{\sl Gap labelling theorems for Schr\"odinger operators},
in From Number Theory to Physics, 538--630, Springer, Berlin, 1992.

\item{[\rLastiii]} Y.~Last,
{\sl Almost everything about the almost Mathieu operator. I},
In: XIth International Congress of Mathematical Physics (Paris, 1994),
pp. 366--372, Cambridge MA: Internat. Press, 1995.

\item{[\rFaKa]} L.D.~Faddeev, R.M.~Kashaev,
{\sl Generalized Bethe ansatz equations for Hofstadter problem},
Commun. Math. Phys. {\bf 169}, 181--191 (1995).

\item{[\rHKW]} Y.~Hatsugai, M.~Kohmoto, Y.-S.~Wu,
{\sl Quantum group, Bethe ansatz equations, and Bloch wave functions in magnetic fields},
Phys. Rev. B {\bf 53}, 9697--9712 (1996).

\item{[\rJito]} S.~Jitomirskaya,
{\sl Metal-insulator transition for the almost Mathieu operator},
Ann. of Math. {\bf 150}, 1159--1175 (1999).

\item{[\rBouJitoii]} J.~Bourgain S.~Jitomirskaya,
{\sl Continuity of the Lyapunov exponent
for quasiperiodic operators with analytic potential},
J. Stat. Phys. {\bf 108}, 1203--1218 (2002).

\item{[\rPuigi]} J.~Puig,
{\sl Cantor spectrum for the almost Mathieu operator},
Commun. Math. Phys. {\bf 244}, 297--309 (2004).

\item{[\rAvKri]} A.~Avila, R.~Krikorian,
{\sl Reducibility or nonuniform hyperbolicity
for quasiperiodic Schr\"odin\-ger cocycles},
Ann. Math. {\bf 164}, 911--940 (2006).

\item{[\rDama]} D.~Damanik,
{\sl The spectrum of the almost Mathieu operator},
Lecture series in the CRC 701 (2008).

\item{[\rGoSch]} M. Goldstein and W. Schlag,
{\sl Fine properties of the integrated density of states
and a quantitative separation property of the Dirichlet eigenvalues},
Geom. Funct. Anal. {\bf 18}, 755--869 (2008).

\item{[\rJiMarx]} S.~Jitomirskaya, C.A.~Marx,
{\sl Analytic quasi-periodic cocycles with singularities
and the Lyapunov exponent of extended Harper's model},
Comm.~Math. Phys. {\bf 316}, 237--267 (2012).

\item{[\rATZ]} C.~Aistleitner, N.~Technau, A.~Zafeiropoulos,
{\sl On the order of magnitude of Sudler products},
Preprint 2020, arXiv:2002.06602.

\item{[\rKochAM]} H.~Koch,
{\sl Golden mean renormalization
for the almost Mathieu operator and related skew products},
J. Math. Phys. {\bf 62}, 1--12 (2021).

\item{[\rKochTrig]} H.~Koch,
{\sl On trigonometric skew-products over irrational circle-rotations},
Discrete Contin. Dynam. Systems A. {\bf 41}, 5455--5471 (2021).

\item{[\rKochUniv]} H.~Koch,
{\sl Asymptotic scaling and universality
for skew products with factors in $\rmSL(2,{\scriptstyle\real})$},
To appear in Erg. Theor. Dyn. Syst.

}

\bye